\documentclass[preprint,12pt]{elsarticle}

 \usepackage{amsthm}
\usepackage{lineno}
\usepackage{latexsym,amscd,amsfonts,enumerate,supertabular}
\usepackage{tabularx}
\usepackage{array}
\usepackage{bigstrut}
\usepackage{graphicx}
\usepackage{epsfig}
\usepackage{amssymb}
\usepackage[intlimits]{amsmath}
\usepackage{algorithm}
\usepackage{multirow}
\usepackage{hhline}
\usepackage{makecell}
\usepackage{color}
\usepackage{xcolor}
\usepackage{float}
\usepackage{graphicx,subfigure}
\usepackage{mathabx}
\usepackage{url}
\usepackage{mathtools, nccmath}
\usepackage{makecell}
\DeclarePairedDelimiter{\nint}\lfloor\rceil

\definecolor{ao(english)}{rgb}{0.0, 0.5, 0.0}

\definecolor{darkmagenta}{rgb}{0.55, 0.0, 0.55}

\newcommand{\ee}{{\rm e}\hspace{1pt}}
\newcommand{\dd}{\hspace{0.5pt}{\rm d}\hspace{0.5pt}}

\newcolumntype{R}[1]{>{\raggedleft\let\newline\\\arraybackslash\hspace{0pt}}m{#1}}
\newcolumntype{L}[1]{>{\raggedright\let\newline\\\arraybackslash\hspace{0pt}}m{#1}}
\newcolumntype{C}[1]{>{\centering\let\newline\\\arraybackslash\hspace{0pt}}m{#1}}

\newcommand{\epiIII}{\mathtt{epi3}}
\newcommand{\exprbIV}{\mathtt{exprb42}}
\newcommand{\pexprbIV}{\mathtt{pexprb43}}
\newcommand{\exprbV}{\mathtt{exprb53}}
\journal{Journal of Computational Physics}

\numberwithin{equation}{section}

\numberwithin{lemma}{section}

\numberwithin{theorem}{section}
\newtheorem{remark}[]{Remark}
\numberwithin{remark}{section}

\numberwithin{example}{section}
\begin{document}
\begin{frontmatter}
\title{{\bf
Further development of the efficient and accurate time integration schemes for meteorological models
}}
 \tnotetext[label1]{This work was supported in part by the
   U.S.~Department of Energy, Office of Science project ``Frameworks,
   Algorithms and Scalable Technologies for Mathematics (FASTMath),''
   under Lawrence Livermore National Laboratory subcontract B626484.}
\author[smuadd]{Vu Thai Luan\corref{cor1}}
 \ead{vluan@smu.edu}
\author[mscadd]{Janusz A. Pudykiewicz}
\ead{Janusz.Pudykiewicz@ec.gc.ca}
\author[smuadd]{Daniel R. Reynolds}
\ead{reynolds@smu.edu}

\cortext[cor1]{Corresponding author}
\address[smuadd]{Department of Mathematics, Southern Methodist University,
  PO Box 750156, Dallas, TX 75275-0156, USA}
\address[mscadd]{Meteorological Service of Canada, Recherche en Pr{\'e}vison
  Num{\'e}rique, 2121 Trans-Canada Highway, Dorval, Que., Canada H9P IJ3}
\begin{abstract}
\small
 In this paper, we investigate the use of higher-order
exponential Rosenbrock time integration methods for the shallow water
equations on the sphere.  This stiff, nonlinear model provides a
`testing ground' for accurate and stable time integration methods in
weather modeling, serving as the focus for exploration of novel
methods for many years.  We therefore identify a candidate set of
three recent exponential Rosenbrock methods of orders four and five
($\exprbIV$, $\pexprbIV$ and $\exprbV$) for use in this model.  Based
on their multi-stage structure, we propose a set of modifications to
the \texttt{phipm/IOM2} algorithm for efficiently calculating the
matrix functions $\varphi_k$.  We then investigate the performance of
these methods on a suite of four challenging test problems, comparing
them against the $\epiIII$ method investigated previously in
\cite{ClancyPudykiewicz13,GaudreaultPudykiewicz16} on these problems.
In all cases, the proposed methods enable accurate solutions at
much longer time-steps than $\epiIII$, proving considerably more
efficient as either the desired solution error decreases, or as the
test problem nonlinearity increases.
\small
\end{abstract}

\end{frontmatter}
\section{Introduction}
\label{sec:introduction}
\def\pop#1#2{ \frac{\partial #1} {\partial #2 } }
\def\dod#1#2{ \frac { d #1} { d #2 }}
\def\oneover #1 #2{ \frac {#1} { #2 }}
\def\bd#1{{\bf #1}}

The idea of predicting the weather by solving fluid equations was developed
at the end of the nineteenth century, long before the appearance of digital computers.
The philosophical basis for the formulation of a forecast problem was deeply
rooted in the tradition of analytical mechanics of Lagrange, Laplace and Jacobi.
The belief in determinism further influenced the formulation of weather prediction
based on the methodology used by Laplace in his analysis of tidal motions.

Despite a strong scientific justification, the meteorological community in the
early 1900s never considered the fluid equations as a serious predictive tool
because the closed form of a solution was unattainable.
This unfavourable reception of a dynamic forecast changed, however, after
Richardson published a book on predicting the weather by a numerical process \cite{Richardson1922}.
This seminal contribution contained the first formulation of algorithms for the
approximate solution of meteorological equations, leading ultimately to the predictions
based on primitive equations which materialized, four decades later, in the mid 1960s.

The widespread use of dynamic models based on the primitive equations led to
new theoretical insights including the discovery of chaos by E.~Lorenz.
The sixties brought to atmospheric sciences the fulfillment of Richardson's dream,
as well as the reformulation of the basic work of H.~Poincar\'e on the stability of
dynamical systems.

Through the entire process of the development of meteorological models,
finding a stable and reasonably efficient time integration scheme has posed a major challenge.
In order to describe the scope of these problems, we
consider the equations of a geophysical rotating fluid
after discretization of the spatial derivatives.
The set of meteorological equations can be cast in the form of an
autonomous dynamical system

\begin{equation}\label{eq:dyn}
\dod u t = F(u), \hskip1.0truecm u(0)=u_0,
\end{equation}
where $u\in \mathbb{R}^n$ is the state vector, $n$ indicates the number
of degrees of freedom, and
$F: \mathbb{R}^n \longrightarrow \mathbb{R}^n$.

The eigenvalues of the Jacobian matrix of $F$ in \eqref{eq:dyn} differ by several orders of magnitude,
reflecting the fact that the primitive meteorological equations govern processes with
time scales ranging from a relatively slow advection to very fast gravity waves.
This property is often described as \emph{stiffness} of the equations.
The difficult task of solving \eqref{eq:dyn} has been a central issue
in the numerical solution of meteorological models over the past 70 years.
At the initial stage of the development of numerical weather prediction models,
the time stepping schemes for integration of the primitive meteorological equations
were explicit and their convergence was subject to severe restrictions imposed by the
Courant Friedrichs Lewy (CFL) stability condition for the gravity waves.

The introduction of semi--implicit methods to meteorological models in the
1970s changed this situation; it became evident that the time step
could be increased sixfold without affecting the overall accuracy
\cite{Robert1969} and 
\cite{Robert1981,Robert1982}.
This significant advancement was based on the observation that
the right-hand side of \eqref{eq:dyn} can be naturally cast in the
following semi--linear form

\begin{equation}\label{eq:dynspl}
\dod u t = L\,u+N(u),
\end{equation}
where $L$ and $N$ are the linear and nonlinear parts, respectively.

The stiffness of \eqref{eq:dyn} comes predominantly from the linear part of $F$.
Consequently, the numerical solution of \eqref{eq:dynspl} in the
semi--implicit approach was performed by approximating the linear term
$L\,u$ implicitly, and the nonlinear part $N(u)$ explicitly.
This methodology offers a compromise between the requirements
of accuracy and efficiency \cite{Robert1986}.  Later, Ascher
  \emph{et.~al.} \cite{Ascher1995} used this idea to develop so-called
  implicit-explicit methods (IMEX) for time-dependent PDEs. IMEX
  methods can be considered as modern, more-accurate, semi-implicit
  methods. They have been widely used in the numerical analysis
  community, and have shown promising results in the shallow water
  context \cite{Kar2006,Bispen2015,WhitakerKar2013}.

Longer time steps afforded by the semi--implicit scheme
liberated significant computing resources which, in turn, were used to
include additional physical processes, including calculating the
parameterization of clouds, boundary layer processes and radiative
transfer.  With further increases of the complexity of models
and an emerging trend to include smaller scale processes,
the semi--implicit scheme was extended to the compressible
Euler equations \cite{Tanguay1990}.  Subsequent numerical experiments
indicated that the semi--implicit scheme could also be applied for the
study of convective scale atmospheric motions \cite {Robert1993}.
The notion of an ``ultimate'' time stepping method for all
meteorological applications was born.

Alternative approaches, based on exponential
time integration method, were proposed in recent years.
After multiplication of \eqref{eq:dynspl} by the integrating
factor, $e^{-Lt}$, we obtain the variation-of-constants formula
\begin{equation}\label{eq:dynst}
u(t)=\ee^{-Lt}u_0 + \int_{0}^{t} \ee^{(t-\tau)L} N(u(\tau)) d\tau.
\end{equation}
Exponential integration schemes based on \eqref{eq:dynspl}-\eqref{eq:dynst}
were first considered in the 1960s \cite{Certaine1960,Lawson1967,Pope63} and later developed by many authors, e.g., Beylkin et al.~\cite{Beylkin1998}, Cox and
Matthews \cite{CoxandMatthews}, and Hochbruck and Ostermann \cite{HO05}. For a recent review on exponential integrators, see \cite{HO10}. Consistent with the fact that the
linear part of the problem is solved exactly, exponential integration schemes
allow the possibility of very good accuracy and realistic
representation of high frequencies, in contrast to semi--implicit schemes.
In meteorological applications, Archibald et al.~\cite{ArchibaldEtAl}
applied the scheme of Beylkin et al.~\cite{Beylkin1998},
based on the assumption of static splitting expressed
by \eqref{eq:dynspl}, to solve the shallow water equations
on the  sphere.

Clancy and Pudykiewicz \cite{ClancyPudykiewicz13}
applied methods based on a \emph{dynamic linearization}
to the shallow water system on an icosahedral geodesic grid \cite{Pudykiewicz2011}.
In this approach, equation \eqref{eq:dyn} is linearized at each time
step with respect to the continuously changing state of the system, $u$.
The $\mathtt{phipm}$ method of Niesen and Wright \cite{NW12}
was used to evaluate exponential functions of the Jacobian operator.
Despite the fact that the method allowed long time steps and was very accurate, it was deemed
to be too expensive to be of practical importance by some users.
This deficiency was addressed in the subsequent paper
by Gaudreault and Pudykiewicz \cite{GaudreaultPudykiewicz16}
who modified the original $\mathtt{phipm}$ algorithm.
They replaced the Arnoldi iteration with the incomplete orthogonalization method
introduced by Saad \cite{Saad1980}.  This strategy followed the idea of using
the incomplete orthogonalization method for the time integration of an
advection--diffusion equation, suggested by Koskela \cite{Koskela2015}.
The efficiency of the exponential scheme described in
\cite{GaudreaultPudykiewicz16} was further enhanced by a new method to
determine the initial size of the Krylov space based on information
from previous time steps.

After introducing optimizations resulting from the application of an incomplete
orthogonalization algorithm, an accurate and efficient numerical integration
of the shallow water equations, with time steps significantly longer than those
in the semi--implicit scheme, became possible.
However, at these larger time steps nonlinear effects became more
prominent, resulting in reduced accuracy due to its simple
approximation of the nonlinear term, $N(u)$.  This fact pointed to the imminent
danger of performing a stable but meaningless integration with a long time step.

The primary purpose of this work is, therefore, to explore time
integration schemes that can retain the positive characteristics of
the methods outlined in \cite{GaudreaultPudykiewicz16}, while offering
a higher-accuracy approximation of the nonlinearity that remains after
dynamic linearization.  One of the most promising candidates to
fulfill theese requirements are the \emph{exponential Rosenbrock
methods} \cite{HO06,HOS09,LO14a,LO16,L17}.  This paper focuses on the
use of such methods for solving the meteorological equations,
including both additional optimizations to the $\mathtt{phipm}$ algorithm (which results in a new routine called $\mathtt{phipm\_simul\_iom2}$)
for exponential Rosenbrock methods, as well as investigations of their
accuracy and efficiency on a range of applicable test problems.
This work is especially relevant in the context of the ongoing
  debate regarding optimal time integration methods for atmospheric
  models \cite{Mengaldo2018}.
The main intention of the search for new methods is not necessarily
the desire to replace well-established methods, but to explore new possibilities.
We hope that they offer certain advantages, such as a better representation of
the phase properties of gravity waves and stiff stability, which will be crucial
when complicated chemistry and micro physics are added to the atmospheric models.

The paper is organized as follows.
The first section provides our rationale for selecting the shallow
water equations for a study of efficient time integration methods.
Following this, we present a concise description of an autonomous system \eqref{eq:dyn}
obtained after the spatial discretization of the governing equations on icosahedral grid.
We also derive an explicit analytical form of the Jacobian operator.
The general formulation of the exponential Rosenbrock schemes, as well
as details regarding their implementation, is covered in Sections 3 and 4.
The extensive discussion of the results obtained with the
new schemes, using the tests reported in \cite{GaudreaultPudykiewicz16}, is included
after the theoretical section.
The main objective of this discussion is to demonstrate that the proposed techniques
are much more accurate and even more efficient than that the algorithms used in the initial tests
reported in \cite{ClancyPudykiewicz13,GaudreaultPudykiewicz16}.
In the last part of the paper, the results are summarized and future work with the
exponential Rosenbrock methods for atmospheric models is outlined.
\section{Shallow Water Model}
\label{sec:shallow_water_model}

\def\mel#1{\pop{#1}{\xi_i} \pop {#1} {\xi_j}}

\def\tth{
\begin{bmatrix}
  -\sin \theta \cos \phi \\
  -\sin \theta \sin \phi \\
    \cos \theta
 \end{bmatrix} }

\def\ttp{
\begin{bmatrix}
 -\cos\theta\sin \phi \\
  \cos\theta\cos \phi \\
  0
 \end{bmatrix} }

The original work with the time integration schemes of Robert \cite{Robert1969}
was performed using a barotropic model of the Earth's atmosphere.
In the years following that seminal contribution, the shallow water equations
became a standard tool for the investigation of prototypes of the meteorological models.
In this study, we will continue to use a barotropic model to develop
more accurate and efficient exponential integration algorithms.

The flow is analyzed on a two-dimensional sphere $S$
described in Cartesian coordinates
by the following parametric equations
\[
  x=a \cos \theta \cos \phi,\, y=a \cos \theta \sin \phi,\, z=a \sin \theta,
\]
where $(\theta,\phi)=(\xi_1,\xi_2)$
are the spherical coordinates,
$\theta$ is the latitude ($\theta\in[-\pi/2,\pi/2]$),
$\phi$ is the longitude ($\phi\in(0,2 \pi])$,
and $a$ is the radius of the sphere.  Specifically, $\theta$ is the
angle measured from the $xy$- plane, and $\phi$ is the angle measured
in the $xy$- plane from the positive $x$ axis, counterclockwise as
viewed from the positive $z$ axis.
At each point of $S$ we can define
the normal vector
\begin{equation}\label{eq:outn}
\bd n=(\cos\theta \cos\phi,\, \cos\theta \sin\phi,\, \sin\theta)^T
\end{equation}
which spans the normal vector space, $TS^{\bot}$.

The tangent space $TS$ is spanned by two basis vectors
\begin{equation}\label{eq:outw}
{\bd e}_1=\pop {\bd n} {\theta}=\tth, \hskip 1cm   {\bd e}_2=\pop {\bd n} \phi=\ttp .
\end{equation}

The Riemannian metric of the sphere with respect to the coordinate system
defined by $({\bd e}_1,{\bd e}_2)$
is given by
\begin{equation}\label{eq:metric}
g_{ij}=\bigg(\mel x + \mel y + \mel z\bigg)= \begin{bmatrix} a^2 & 0 \\ 0&a^2\cos^2\theta \end{bmatrix},
\end{equation}
and the inverse of $(g_{ij})$ is denoted by $(g^{ij})$, i.e.~$g_{ij}\,g^{jk}=\delta^k_i$.

The equations describing flow of a thin layer of fluid on $S$ can be cast
in the form
\begin{align}
  \label{eq:momentumstc}
  \pop {\bd u} t = -({\rm Curl}_n \bd u + &f \bd n) \times \bd u  -{\rm grad} \bigg(\oneover {|\bd u|^2} 2 +g(h+h_s)\bigg), \\
  \label{eq:continuity}
  &\pop {h} t + {\rm div} ({h} \bd u) = 0,
\end{align}
where $\bd u$ is the smooth velocity field on $S$ with values in $TS$,
$h$ is the smooth scalar field on $S$ describing the thickness of the fluid layer,
$h_s$ is the height of the surface level ($h_s \ll a$),
$g$ is the gravitational acceleration,
$f=2 \Omega \sin\theta$ is the Coriolis parameter,
$\Omega$ is the angular velocity of the rotation,
and
$\times$ denotes the vector product in $\mathbb{R}^3$.

After expressing the shallow water equations
\eqref{eq:momentumstc}-\eqref{eq:continuity} in terms of the Cartesian
components they are discretized on the spherical geodesic grid
using the finite volume method.  A detailed description of the discretization,
including the method used to optimize the mesh, is presented in \cite{Pudykiewicz2011}.
The system of the ordinary differential equations obtained after integration
over the control volumes, and applying the discretization formulae for
approximation of the differential operators, can be written in the following
compact form
\begin{align}
  \notag
  \dod {u_x} t &= - \eta \, \cdot P_x  - G_x \cdot {\cal E} \\
  \label{eq:sweqsinv}
  \dod {u_y} t &= - \eta \, \cdot P_y  - G_y \cdot {\cal E} \\
  \notag
  \dod {u_z} t &= - \eta \, \cdot P_z  - G_z \cdot {\cal E} \\
  \notag
  \dod {h} t &= -D_x\cdot(u_x\,h) -D_y\cdot(u_y\,h) -D_z\cdot(u_z\,h),
\end{align}
where $u_x$, $u_y$, $u_z$, ${\cal E}=(u_x^2+u_y^2+u_z^2) / 2 +g(h+h_s)$, and $h$
are the column arrays with Cartesian components of the
velocity ${\bf u}$, total energy and the height field, respectively
(all quantities are volume averages over the control volume),
$G_x$, $G_y$, $G_z$ are sparse matrices used to evaluate the
Cartesian components of the gradient on the sphere, while $D_x$,
$D_y$, and  $D_z$ are sparse matrices used to evaluate the divergence
${\rm div} ({\bf u})=D_x \cdot u_x+D_y \cdot u_y+D_z \cdot u_z$, and
$\eta$ is a column array containing the control volume average
values of the absolute vorticity,
\begin{equation}\label{eq:vort}
\eta= {\rm Curl}_n \bd u  +f=V_x \cdot u_x + V_y \cdot u_y + V_z \cdot u_z +f.
\end{equation}
Sparse arrays $V_x$, $V_y$, and $V_z$ are used to calculate the
vertical component of the vorticity.  For further details consult \cite{Pudykiewicz2011}.
The arrays $P_x$, $P_y$, and $P_z$ contain the Cartesian components of
the vector product of the velocity and the surface normal, and have the form
\begin{align*}
  P_x &= n_y u_z - n_z u_y, \\
  P_y &= n_z u_x - n_x u_z, \\
  P_z &= n_x u_y - n_y u_x,
\end{align*}
where $n_x,\,n_y$, and $n_z$ are column arrays containing the $x,\,y$,
and $z$ components of the normal vector $\bd n$ evaluated at the
control volume centers.

A dissipation term of the form ${D_f}\cdot\psi$ is added for each
prognostic variable $\psi$ in \eqref{eq:sweqsinv}.
The dissipation operator is the same as that used
in the experiments reported in \cite{ClancyPudykiewicz13},
\begin{equation}\label{eq:swdis}
{D_f}= - \nu {{\cal L}}^2,
\end{equation}
where ${\cal L}$ is a sparse matrix representing the Laplace operator.
The dissipation coefficient $\nu$ in \eqref{eq:swdis} is given as
\begin{equation}\label{eq:niu}
\nu=\gamma_{\rm h} \frac{ \left(\overline{\Delta x}\right)^{n_{\gamma}}}{\Delta t_e},
\end{equation}
where $\gamma_{\rm h}$ is the coefficient of proportionality, $n_{\gamma}=4$,
$\Delta t_e=240 s$,
\[
  \overline{\Delta x}=\sqrt{{{4 \pi a^2}/{ N_{\rm g}}}}
\]
is the average separation of the node points
and $N_{\rm g}$ is the number of the nodes in the geodesic grid
($N_{\rm g}=10 \times 2^{2 l} +2 $ with $l$ indicating the grid number).

The system of equations \eqref{eq:sweqsinv} can be rewritten in the
more compact form of the autonomous system \eqref{eq:dyn}, with a
state vector $u$ containing the components of the velocity and the
height field as
\begin{equation}\label{eq:dynu}
u=(u_x^T \hskip 0.2truecm  u_y^T \hskip 0.2truecm u_z^T \hskip 0.2truecm h^T)^{^T}.
\end{equation}
The dimension of the state vector $u$ is thus $4\,N_g$.

The Jacobian $J$ of the autonomous system \eqref{eq:sweqsinv}
is a sparse matrix of size $(4 N_g) \times (4 N_g)$, and is
given by the formula
\begin{equation}\label{eq:dynj}
J=\frac{\partial F} {\partial u}.
\end{equation}

This Jacobian can be written in the symbolic form
\begin{equation}\label{eq:dynjs}
J=J_v+J_r+J_m,
\end{equation}
where the block array $J_v$ represents relative vorticity
and the block array $J_r$ accounts for rotation and dissipation.
The third block array  $J_m$ represents
coupling between the mass field and the velocity field.
These block arrays can be written as
\begin{equation}\label{eq:dynja}
J_v= -\begin{bmatrix} \{P_x^i\,\,V_x^{ij}\}    & \{P_x^i\,\,V_y^{ij}\}  & \{P_x^i\,\,V_z^{ij}\} & \{0\} \\
                                                                                                                                                         \\
                          \{P_y^i\,\,V_x^{ij}\}    & \{P_y^i\,\,V_y^{ij}\}  & \{P_y^i\,\,V_z^{ij}\} & \{0\}  \\
                                                                                                                                                          \\
                          \{P_z^i\,\,V_x^{ij}\}    & \{P_z^i\,\,V_y^{ij}\}  & \{P_z^i\,\,V_z^{ij}\} & \{0\}  \\
                                                                                                                                                          \\
                          \{0\}    &  \{0\}   &  \{0\} &  \{0\} \end{bmatrix},
\end{equation}
\begin{equation}\label{eq:dynjb}
J_r= \begin{bmatrix} \{D_f^{ij}\}    & \{{\eta}^i\,\,n_z^j\,\,{\delta}^{ij}\}  & -\{{\eta}^i\,\,n_y^j\,\,{\delta}^{ij}\} & \{0\} \\
                                                                                                                                                                                    \\
                           \{{\eta}^i\,\,n_z^j\,\,{\delta}^{ij}\}    & \{D_f^{ij}\}  & \hskip 5pt -\{\eta^i\,\,n_x^j\,\,{\delta}^{ij}\} & \{0\}  \\
                                                                                                                                                                                    \\
                          \{{\eta}^i\,\,n_y^j\,\,{\delta}^{ij}\}    & -\{{\eta}^i\,\,n_x^j\,\,{\delta}^{ij}\} & \{D_f^{ij}\} & \{0\}  \\
                                                                                                                                                                                   \\
                          \{0\}    &  \{0\}   &  \{0\} &  \{D_f^{ij}\} \end{bmatrix},
\end{equation}
and
\begin{equation}\label{eq:dynjc}
J_m= -\begin{bmatrix} \{G_x^{ij}\,\,u_x^j\}    & \{G_x^{ij}\,\,u_y^j\}   & \{G_x^{ij}\,\,u_z^j\} & \{g\,G_x^{ij}\} \\
                                                                                                                                                                                    \\
                         \{G_y^{ij}\,\,u_x^j\}    & \{G_y^{ij}\,\,u_y^j\}   & \{G_y^{ij}\,\,u_z^j\} & \{g\,G_y^{ij}\} \\
                                                                                                                                                                                    \\
                         \{G_z^{ij}\,\,u_x^j\}    & \{G_z^{ij}\,\,u_y^j\}   & \{G_z^{ij}\,\,u_z^j\} & \{g\,G_z^{ij}\}  \\
                                                                                                                                                                                    \\
                         \{D_x^{ij}\,\,h^j\}    &  \{D_y^{ij}\,\,h^j\}    & \{D_z^{ij}\,\,h^j\}  & \{D_x^{ij}\,\,u_x^j+D_y^{ij}\,\,u_y^j+
                           D_z^{ij}\,\,u_z^j\}  \end{bmatrix},
\end{equation}
where the wave brackets denote matrices (e.g., the expression $
\{P_x^i\,\,V_x^{ij}\} $ represents the matrix with $ij$-th element
equal to $P_x^i\,\,V_x^{ij}$).
\section{Exponential Rosenbrock methods}
\label{sec:exponential_methods}
In this section, we first recall the idea behind exponential
Rosenbrock methods and display the stiff order conditions for methods
of order up to 5.  We then select a set of promising schemes for our
numerical experiments of the shallow water equations on the
sphere. These include a fourth-order two-stage scheme, a fourth-order
parallel stages scheme,  and a fifth-order three-stage scheme. Our
main references in this section are \cite{HO06,HOS09,LO14a,LO16,L17}.

\subsection{General motivation and ideas}
As discussed in the introduction, when integrating the full nonlinear
stiff system (1), classical methods like explicit Runge--Kutta are subject to the CFL condition, resulting in
unrealistically small time steps. This is mainly because for such
systems the Jacobian of the forcing term $F(u)$ often has a large norm
or is even an unbounded operator (causing the stiffness of the
system). To overcome this, implicit methods are often used. These
standard methods, however, require the solution of nonlinear systems
of equations at each step. As this stiffness increases, such methods
may require increased computational effort.  Alternatively, we
consider a class of integrators that can handle the stiffness of this
full nonlinear system in an explicit and very accurate way. The idea is first to replace
the full nonlinear system (1) by a sequence of semilinear problems
(similarly to the idea of deriving Rosenbrock-type methods, see
\cite[Chap. IV.7]{HW96}). This can be done by linearizing the flow in
each time step around the numerical solution $u_n$ (due to
\cite{Pope63}), leading to
\begin{equation} \label{eq3.1}
  u'(t)=F(u(t))=J_n u(t)+N_n(u(t)),
\end{equation}
where
\begin{equation} \label{eq3.2}
  J_n =\frac{\partial F} {\partial u}(u_n) \quad\text{and}\quad N_n(u)=F(u)-J_n u
\end{equation}
are the Jacobian and the nonlinear remainder, respectively.
One can then apply exponential Runge--Kutta methods (see \cite{HO05})
to the semilinear system \eqref{eq3.1}; these handle the stiffness by
solving the linear part $J_n u$ exactly and integrating the
nonlinearity $N_n (u)$ (which is much smaller than the original
$F(u)$) in an explicit manner.  In this regard, the exact solution at
time $t_{n+1}=t_n+\Delta t $ of \eqref{eq3.1} may be represented using the
variation-of-constants formula \eqref{eq:dynst},
\begin{equation} \label{eq3.3}
u(t_{n+1})=\ee^{\Delta t J_n}u(t_n) +\int_{0}^{\Delta t} \ee^{(\Delta t-\tau)J_n} N_n(u(t_n+\tau)) \dd \tau.
\end{equation}
This structure provides a recipe for constructing these integrators:
the linear part can be integrated exactly by evaluating the action of
the matrix exponential $\ee^{\Delta t J_n}$ on the vector $u(t_n)$, and the integral
involving the new nonlinearity $N_n (u)$ can be approximated by some
quadrature.
Details on the quadrature method used in this work are provided in Section 3.2.
We note that due to the Rosenbrock structure, integration
of this nonlinear term can leverage the fact that
$\frac{\partial N_n} {\partial u}(u_n)=\frac{\partial F} {\partial u}(u_n)-J_n=0$.  This overall procedure results in the
so-called \emph{exponential Rosenbrock methods}, see
\cite{HO06,HOS09}.  It is also worth mentioning that other classes of
exponential integrators were constructed based on \eqref{eq3.3} as
well, see \cite{HO10}.

\begin{remark}
\label{remark3.1}
\textnormal{
In many applications, the forcing term $F(u)$ naturally has the
semilinear form $F(u(t))=Lu(t) +N(u(t))$. Thus stiff semilinear
problems of the form  \eqref{eq:dynspl} can be considered as a
fixed linearization problem ($J_n=L$). Again, one can approximate its
solution by \eqref{eq3.3} as described above, which results in the
so-called \emph{explicit exponential Runge--Kutta methods} (see
\cite{HO05}). While these methods do not suffer from the CFL condition
for the linear part, the stepsize $\Delta t $ is still limited by the CFL
condition due to the nonlinear part. Hence when $N(u)$ is large, the
product of $\Delta t $ and the Lipschitz constant of $N(u)$ should be
sufficiently small to ensure linear stability. In this context, dynamic
linearization along the numerical solution offers a potential
advantage.  This is because in the splitting of the right-hand side,
\begin{equation} \label{eqRemark3.1}
u'(t)=Lu(t) +N(u(t))=J_n u(t)+N_n(u(t)),
\end{equation}
dynamic linearization constructs $J_n = A+N'(u_n)$, leading to a much
smaller nonlinearity $N_n(u)$ that has reduced Lipschitz constant
since $N'_n(u_n)=0$.  We further note that the new linear part $J_n u
= \left(A+N'(u_n)\right)u$ can again be solved exactly.  We therefore
anticipate that exponential Rosenbrock methods may enable use of even
larger time steps than standard explicit exponential Runge--Kutta
methods.  Moreover, another advantage of exponential Rosenbrock
methods is that $N'_n(u_n)=0$, which considerably simplifies the
order conditions, and in turn the derivation of higher-order methods.
}
\end{remark}

To illustrate this idea, we present a simple derivation of a
second-order scheme.

\subsection{A second-order scheme and general schemes}
The integral in \eqref{eq3.3} is approximated as follows.
First, we expand $u(t_n+\tau)$ in a Taylor series as
$u(t_n+\tau)=u(t_n)+\tau u'(t_n)+ r_{u}(\tau)$ with the remainder term
\[
r_{u}(\tau)=\int_{t_n}^{t_n + \tau} u''(\xi) (t_n + \tau-\xi) \dd \xi.
\]
 We then
insert this into $N_n(u(t_n+\tau))$, and perform another Taylor
expansion around $u(t_n)$ (leveraging the fact that $N'_n(u(t_n))=0$),
to obtain
\begin{equation} \label{eq3.4}
N_n (u(t_n+\tau))=N_n (u(t_n)) +r_{N}(\tau)
\end{equation}
with the remainder term
\[
r_{N}(\tau)=\int_{u(t_n)}^{u(t_n + \tau)} N''_n (\kappa) \big(u(t_n)+\tau u'(t_n)+r_{u}(\tau)-\kappa \big)\dd\kappa.
\]
Inserting \eqref{eq3.4} into  \eqref{eq3.3} and denoting $\varphi_1 (\Delta t J_n)=\frac{1}{\Delta t} \int_{0}^{\Delta t} \ee^{(\Delta t -\tau) J_n} \dd \tau$ gives
\begin{equation}\label{eq3.4a}
u(t_{n+1})=\ee^{\Delta t J_n} u(t_n) +\Delta t \varphi_1 (\Delta t J_n) N_n(u(t_n))  +R_n(\Delta t),
\end{equation}
where $R_n(\Delta t)$ is the remainder term, which is given by
\begin{equation}\label{eq:remainder}
R_n(\Delta t)=\int_{0}^{\Delta t} \ee^{(\Delta t-\tau)J_n} r_{N}(\tau)\dd\tau.
\end{equation}
Under standard regularity assumptions on $u(t)$, $N_n(u)$, and the Jacobian $J=\frac{\partial F} {\partial u}$ (mentioned at the
end of this subsection), one can show that $\|r_{u}(\tau)\| \leqslant C \tau^2$
and thus $\|r_{N}(\tau)\| \le C \tau^2$ (e.g., using the Mean Value Theorem for
Integrals). Therefore, it is clear that $\|R_n(\Delta t)\| \leqslant C \Delta t^3$
(where $C$ is a generic constant that may have different values at different
occurrences). We will represent such remainder terms using the Landau
notation, i.e., $R_n(\Delta t)= \mathcal{O}(\Delta t^3)$.  As seen from
\eqref{eq:remainder} the constant $C$ behind $\mathcal{O}$ only depends on
values that are uniformly bounded by the regularity assumptions on $u(t)$ and
$N_n(u)$, but is independent of $\|J_n\|$.

Neglecting the local error term $R_n(\Delta t)= \mathcal{O}(\Delta t^3)$ in  \eqref{eq3.4a} leads to a
second-order scheme, which can be reformulated as
\begin{equation} \label{eq3.5}
u_{n+1}=u_n + \Delta t \varphi_1(\Delta t J_n)F(u_n)
\end{equation}
by replacing $N_n(u(t_n))$ with \eqref{eq3.2} and calculating $\varphi_1(z)=(e^z-1)/z$.
Since as $J_n\rightarrow \mathbf{0}$, \eqref{eq3.5} approaches the
explicit Euler method, it is so-called the \emph{exponential Rosenbrock-Euler method}.
We note that this method has been derived before, e.g.,
\cite{HO06,Hochbruckexp4} that use a different construction.  In
\cite{ClancyPudykiewicz13,Tok06} this method is named $\mathtt{EPI2}$.
The present derivation, however, shows directly that this scheme has
a local error (consistency) of order 3 and thus is a second-order stiffly accurate
method.

To derive higher-order schemes, one must build up higher-order
approximations of the integral in \eqref{eq3.3}, i.e., with a
remainder term of $\mathcal{O}(\Delta t^q)$ ($q\geq 4$). For example, one can
approximate $N_n (u(t_n+\tau))$ by its Taylor expansion of
higher-order and plug it into \eqref{eq3.3}, which introduces the
family of $\varphi$ functions (similar to $\varphi_1$)
\begin{equation} \label{eq3.6}
\varphi _{k}(\Delta t Z)=\frac{1}{\Delta t^k}\int_{0}^{\Delta t} \ee^{(\Delta t-\tau )Z} \tau^{k-1} \dd \tau , \quad k\geq 1,
\end{equation}
which are bounded (see e.g. \cite{HO10}) and satisfy the recursion relation
\begin{equation} \label{eq3.7}
\varphi _{k+1}(z)=\frac{\varphi _{k}(z)-\frac{1}{k!}}{z}, \quad k\geq 0, \ \text{where} \ \varphi_0(z)=\ee^{z}.
\end{equation}
This approach, however, requires the computation of higher derivatives
of $N_n (u)$ which can suffer from instabilities, particularly for
schemes of order higher than three, see \cite{Koskela13}. Therefore, a
preferable solution is to approximate the integral in \eqref{eq3.3} by
using a higher-order quadrature rule, yielding
\begin{equation} \label{eq3.8}
u(t_{n+1})\approx \ee^{\Delta t J_n}u(t_n) +  \Delta t \sum_{i=1}^{s} b_i(\Delta t J_n)N_n( u(t_n+c_i \Delta t))
\end{equation}
with nodes $c_i$ in $[0,1]$  and weights which are matrix functions of
$\Delta t J_n$, denoted by $b_i(\Delta t J_n)$ (similarly to the construction of
classical Runge--Kutta methods). This introduces a set of unknown
intermediate values, $u(t_n+c_i \Delta t)$, which can again be approximated
by using \eqref{eq3.3} (with $c_i \Delta t$ in place of $\Delta t $). To avoid
generating new unknowns, another quadrature rule with the same nodes
$c_j$ ($1\leq j \leq i-1$ for explicit schemes) and new weights
$a_{ij}(\Delta t J_n)$ is used,
\begin{equation} \label{eq3.9}
u(t_n+c_i \Delta t)\approx \ee^{c_i \Delta t J_n}u(t_n) + \Delta t \sum_{j=1}^{i-1} a_{ij}(\Delta t J_n)N_n(u(t_n+c_j \Delta t)).
\end{equation}
As was done for \eqref{eq3.5}, we may reformulate \eqref{eq3.8} and
\eqref{eq3.9} to obtain the general format of $s$-stage explicit
exponential Rosenbrock methods
\begin{subequations} \label{eq3.10}
\begin{align}
U_{ni}&= u_n + c_i \Delta t \varphi _{1} ( c_i \Delta t J_n)F(u_n) + \Delta t \sum_{j=2}^{i-1}a_{ij}(\Delta t J_n) D_{nj}, \label{eq3.10a} \\
u_{n+1}& = u_n + h\varphi _{1} (\Delta t J_n)F(u_n) +  \Delta t \sum_{i=2}^{s}b_{i}(\Delta t J_n) D_{ni}  \label{eq3.10b}
\end{align}
with
\begin{equation} \label{eqDni}
  D_{ni}= N_n ( U_{ni})- N_n(u_n ), \ \  2\leq i\leq s
\end{equation}
\end{subequations}
where $u_n \approx u(t_n)$ and $U_{ni}\approx u(t_n +c_i \Delta t_n)$ (one can define $U_{n1}=u_n$ and consequently $c_1=0$ since they do not enter the scheme due to the reformulation, see \cite{LO16}).
Similarly to \eqref{eq3.4}, since $N_n (U_{ni}) -
N_n(u_n)=\mathcal{O}(\Delta t^2)$ then we have $D_{ni}=\mathcal{O}(\Delta t^2)$ as
well. Thus as seen from \eqref{eq3.10}, the general methods are small
perturbations of the exponential Rosenbrock-Euler method
\eqref{eq3.5}.
By construction, the weights $a_{ij}(\Delta t J_n)$ and $b_i(\Delta t J_n)$ are
usually linear combinations of the $\varphi$ functions, $\varphi _{k}
(c_i \Delta t J_n)$ and $\varphi_{k} (\Delta t J_n)$, respectively.  Therefore, they
are also uniformly bounded independently of $\|J_n\|$ (i.e.~the
stiffness). This is a very important feature in comparison to
classical exponential schemes (e.g., \cite{CoxandMatthews,Tok06}), where the
matrix functions are expanded using \emph{classical} Taylor series
expansions (involving powers of the Jacobian, e.g.~$(J_n)^k$), that
are only valid for nonstiff problems with small $\|J_n\|$.

Another significant advantage of exponential Rosenbrock methods
\eqref{eq3.10} is that they are fully explicit, and do not require the
solution of linear or nonlinear systems of equations.

In general, the convergence of exponential Rosenbrock
methods is analyzed for stiff problems where the Jacobian $J=\frac{\partial F} {\partial u}$ generates a strongly continuous semi-group in some Banach space $X$ (so  $\|\ee^{tJ}\|\leq C$ holds uniformly) and furthermore one needs some regularity assumptions on the solution $u(t)$
(sufficiently smooth) and $N_n(u)$ (sufficiently Fr\'echet
differentiable in a neigborhood of the solution), with uniformly
bounded derivatives. For more details, see \cite{HOS09,LO14a}.

\subsection{Selected exponential schemes for numerical experiments}
\label{sec:selected_exponential_schemes}
As seen above, in order to derive exponential Rosenbrock schemes for
stiff problems, one must determine the coefficients $a_{ij}(\Delta t J_n)$ and
$b_i(\Delta t J_n)$ that guarantee approximation of the nonlinear term to
a desired accuracy, even in the presence of large $\|J_n\|$.
Therefore, a very careful local error analysis must be performed to
make sure that the error terms do not contain powers of $J_n$.
Recently, Luan and Ostermann \cite{LO13} derived a new stiff order
conditions theory and performed a convergence analysis for
methods of arbitrary order.  In Table~\ref{tb3.1} we display the
required 4 conditions for deriving methods up to order 5, which can be
also found in \cite{LO14a}. Note that, for exponential Runge--Kutta
methods (applied to a fixed linearization), one needs 16 order
conditions for methods up to order 5 (see \cite{LO12b}). This confirms
the observation in Remark~\ref{remark3.1} regarding the advantage of
using exponential Rosenbrock methods, based on dynamic linearization
of the ODE right-hand side.
\vspace{-5pt}
\renewcommand{\arraystretch}{1.4}%
\begin{table}[H]
\caption{Stiff order conditions for exponential Rosenbrock methods up to order 5. Here $Z, K$ denote arbitrary square matrices and $\psi_{3,i}(z)= \sum_{k=2}^{i-1}a_{ik}(z)\frac{c^2_k}{2!}-c^{3}_i \varphi_{3} ( c_i z)$. }
\begin{center}
\begin{tabular}{ |c|c|c| }
\hline
No. & Order condition & Order \\
\hline
1&$\sum_{i=2}^{s} b_i (Z)c^2_i=2\varphi_3 (Z) $&3 \\
 \hline
2&$\sum_{i=2}^{s} b_i (Z)c^3_i=6\varphi_4 (Z) $&4 \\
\hline
3&$\sum_{i=2}^{s} b_i (Z)c^4_i=24\varphi_5 (Z) $&5 \\
4&$\sum_{i=2}^{s} b_i (Z)c_i  K \psi _{3,i}(Z)=0 $&5 \\
\hline
\end{tabular}
\end{center}
\label{tb3.1}
\end{table}%
\vspace{-10pt}
With these stiff order conditions in hand, one can easily derive
numerous methods of order up to 5. For practical implementation,
however, one must optimize the coefficients and stages for improved
accuracy and computational efficiency.  Guided by this, we select the
following three representative schemes for our applications.

First, we consider a fourth-order scheme satisfying the stiff order
conditions,
named $\mathtt{exprb42}$ in \cite{L17}: 
\begin{subequations} \label{eq3.11}
\begin{align}
 U_{n2}&= u_n + \tfrac{3}{4} \Delta t \varphi _{1} ( \tfrac{3}{4}\Delta t J_n)F(u_n), \label{eq3.11a} \\
u_{n+1}& = u_n + \Delta t \varphi _{1} (\Delta t J_n)F(u_n) \\
       &+ \Delta t \tfrac{32}{9}\varphi_3 (\Delta t J_n) (N_n (U_{n2})- N_n(u_n )) \label{eq3.11b}.
\end{align}
\end{subequations}
We note that this uses only two stages, and is therefore considered as a superconvergent scheme.

Second, we consider a fourth-order 3-stage scheme satisfying the stiff order conditions, named $\mathtt{pexprb43}$ in
\cite{LO16}:
\begin{subequations} \label{eq3.12}
\begin{align}
U_{n2}&= u_n + \tfrac{1}{2} \Delta t \varphi _{1} ( \tfrac{1}{2}\Delta t J_n)F(u_n), \\
U_{n3}&= u_n +  \Delta t \varphi _{1} (\Delta t J_n)F(u_n), \\
u_{n+1}& = u_n + \Delta t \varphi _{1} (\Delta t J_n)F(u_n) + \Delta t \varphi_3 (\Delta t J_n) (16 D_{n2}-2 D_{n3})\\
          &+  \Delta t \varphi_4 (\Delta t J_n) (-48 D_{n2}+12 D_{n3}).
\end{align}
\end{subequations}
We note that since $U_{n3}$ does not depend on $U_{n2}$ these two
stages may be computed in parallel (hence the ``p'' preceding the name). \\

Third, we consider a fifth-order scheme satisfying the stiff order
conditions, named $\mathtt{exprb53}$ in \cite{LO14a}: 
\begin{subequations}  \label{eq3.13}
\begin{align}
 U_{n2}&= u_n + \tfrac{1}{2} \Delta t \varphi _{1} ( \tfrac{1}{2}\Delta t J_n)F(u_n), \\
  U_{n3}&= u_n + \tfrac{9}{10} \Delta t \varphi _{1} ( \tfrac{9}{10}\Delta t J_n)F(u_n)\\
  &+ \Delta t \tfrac{27}{25}\varphi _{3} ( \tfrac{1}{2}\Delta t J_n) D_{n2} +\Delta t \tfrac{729}{125}\varphi _{3} ( \tfrac{9}{10}\Delta t J_n) D_{n2}, \\
u_{n+1}& = u_n + \Delta t \varphi _{1} (\Delta t J_n)F(u_n) + \Delta t \varphi_3 (\Delta t J_n) (18D_{n2}-\tfrac{250}{81}D_{n3})\\
&+  \Delta t \varphi_4 (\Delta t J_n) (-60D_{n2}+\tfrac{500}{27}D_{n3}),
\end{align}
\end{subequations}
We note that this is also superconvergent, as it uses only three
stages.  \\

We further note that the vectors $D_{n2}$ and $D_{n3}$ in the methods
$\mathtt{pexprb43}$ and $\mathtt{exprb53}$ are given by \eqref{eqDni},
i.e.,
 $D_{n2}=N_n (U_{n2})- N_n(u_n ), \ \ D_{n3}=N_n (U_{n3})- N_n(u_n )$.

For comparison purposes, we also include the third-order exponential
multistep scheme named $\mathtt{epi3}$, proposed in \cite{Tok06}: 
 \begin{equation}  \label{eq_epi3}
u_{n+1} = u_n + \Delta t \varphi _{1} (\Delta t J_n)F(u_n)+\Delta t \tfrac{2}{3}\varphi _{2} (\Delta t J_n)R_{n-1},
\end{equation}
where $R_{n-1}=F(u_{n-1})-F(u_n)-J_{n} (u_{n-1}-u_n)$.  This scheme is
a classical method as it does not satisfy the \emph{stiff} order
conditions for exponential multistep methods of order 3 (see \cite{HO11}), i.e., the
error terms contain powers of $J_n$.  Given the initial value $u_0$,
this method requires the computation of the additional starting value
$u_1$.  For the numerical results in Section
\ref{sec:numerical_experiments}, we compute this additional initial
condition using a single step of the exponential Rosenbrock-Euler
method \eqref{eq3.5}; as the local error of this method matches the
$\mathcal O(\Delta t^3)$ global error of \eqref{eq_epi3}, it provides a
suitable method for generation of this missing initial condition.  We
finally note that this method was previously used in
\cite{ClancyPudykiewicz13,GaudreaultPudykiewicz16} for the simulation
of the shallow water model on the sphere, and so inclusion of this
method in our tests allows direct comparison with the methods tested
in those earlier works.

\section{Implementation of exponential integrators}
\label{sec:implement}
In this section, we briefly review some state-of-the-art algorithms
for the implementation of exponential integrators. We then discuss our
modifications to one of these algorithms to enhance efficiency for
our proposed exponential Rosenbrock methods.  Finally, we discuss
specific details in using our new algorithm for the three selected
exponential Rosenbrock methods $\mathtt{exprb42}$,
$\mathtt{pexprb43}$, and $\mathtt{exprb53}$.

\subsection{The state-of-the-art algorithms}
The implementation of exponential integrators requires computing the
action of matrix functions $\varphi_k (\Delta t J_n)$ on vectors $v_k$.  In
1997, Hochbruck and Lubich \cite{HL97} analyzed Krylov subspace
methods for efficiently computing the action of a matrix exponential
(with large norm) on some vector.  Since then, there have been dramatic
advances in constructing more efficient techniques, even for the
evaluation of a linear combinations of $\varphi$-functions acting on
sets of vectors $v_0,\ldots,v_p$,
\begin{equation} \label{eq3.14}
\varphi_0 (A)v_0 + \varphi_1 (A)v_1 +\varphi_2 (A)v_2+\cdots+\varphi_p(A)v_p,
\end{equation}
which is crucial within calculations of each stage in exponential schemes.
In particular, we highlight a number of state-of-the-art algorithms.
$\mathtt{expmv}$, proposed by Al-Mohy and Higham in \cite{AH11}, uses
a truncated standard Taylor series expansion. $\mathtt{phipm}$ was
proposed by Niessen and Wright in \cite{NW12}, and uses adaptive
Krylov subspace methods.  Finally, $\mathtt{expleja}$ was proposed by
Caliari \emph{et al.}~in \cite{CKOS16}, and uses Leja interpolation.
Among these, it turns out that $\mathtt{phipm}$ offers an
advantage in terms of computational time.  This algorithm is developed
based on an adaptive time-stepping method to evaluate \eqref{eq3.14}
using only one matrix function, $\varphi_p (\tau_k A)w_p$, in each
substep, where $w_p$ is a vector resulting from the time-stepping
method, and $\tau_k\le 1$ is a substep size.  This can be carried out
in a lower dimensional Krylov subspace (with one Krylov projection
needed).  When computing this matrix function, standard Krylov
subspace projection methods are employed, i.e., using the Gram-Schmidt
orthogonalization procedure (the Arnoldi iteration).  For improved
accuracy and efficiency, the dimension of these Krylov subspaces and
the number of substeps $\tau_k$ are chosen adaptively.  More recently,
the $\mathtt{phipm}$ routine was modified by Gaudreault and
Pudykiewicz in \cite{GaudreaultPudykiewicz16} to use the incomplete
orthogonalization method (IOM) within the Arnoldi iteration.  They
additionally adjusted the two crucial initial parameters for starting
the Krylov adaptivity in order to further reduce the computational
time. Their result is a new routine named $\mathtt{phipm/IOM2}$.  It
is shown in \cite{GaudreaultPudykiewicz16} that this algorithm offers
a significant computational advantage compared to $\mathtt{phipm}$ for
integrating the shallow water equations on the sphere.

\subsection{The time-stepping technique for computing a linear combination of $\varphi$-functions}
In order to contextualize the further improvements that we have made
to $\mathtt{phipm/IOM2}$ described in Section
\ref{sec:implementation}, we describe the idea of the adaptive
time-stepping method (see \cite{NW12}) for evaluating the linear
combination \eqref{eq3.14} efficiently.

The linear combination \eqref{eq3.14} is precisely equal to the solution
$y(1)$ of the ODE
\begin{equation} \label{eq3.15}
y'(t)=Ay(t)+v_1+t v_2+\cdots+\frac{t^{p-1}}{(p-1)!}v_p, \ y(0)=v_0.
\end{equation}
The time-stepping technique leverages this equivalence to approximate
$y(1)$ by discretizing $[0,1]$ into subintervals $0=t_0<t_1<\cdots<t_K=1$, having
widths $\tau_k=t_{k+1}-t_k$.  Then by exploiting the relation between $y(t_{k+1})$ and $y(t_k)$ (see \cite[Lemma 2.1]{NW12}), we have
\begin{equation} \label{eq3.16}
y(t_{k+1})=\varphi_0 (\tau_k A)y(t_k)+\sum_{i=1}^{p}\tau^i_k \varphi_i (\tau_k A) \sum_{j=0}^{p-i} \frac{t^j_k}{j!}v_{i+j}.
\end{equation}
With the help of the recursion relation \eqref{eq3.7}, \eqref{eq3.16}
may be simplified to
\begin{equation} \label{eq3.17}
y(t_{k+1})=\tau^p_k \varphi_p (\tau_k A)w_p+ \sum_{j=0}^{p-i} \frac{\tau^j_k}{j!}w_j,
\end{equation}
where the vectors $w_j$ satisfy another recurrence relation,
\begin{equation} \label{eq3.18}
w_0=y(t_k), \ w_j=Aw_{j-1}+\sum_{\ell=0}^{p-j} \frac{t^{\ell}_k}{\ell!}v_{j+\ell}, \ j=1,\ldots,p.
\end{equation}
This interesting result implies that evaluating the linear combination
\eqref{eq3.14} (i.e.~$y(t_K)=y(1)$), which typically consists of
$(p+1)$ matrix-vector multiplications, can be accomplished instead by
using only one matrix function $\varphi_p (\tau_k A)w_p$ in each
substep. Since $0<\tau_k\le 1$, this task can be carried out in a
lower dimensional Krylov subspace, and only one Krylov projection is
needed.  Therefore, this approach may reduce the computational cost
significantly in comparison with a standard Krylov algorithm.  More
specifically, it is known that constructing a Krylov subspace of
dimension $m$ for approximating $\varphi_p (A)v$ requires a total
computational cost of $\mathcal{O}(m^2)$.  The time-stepping method
\eqref{eq3.17}-\eqref{eq3.18}, replaces this by $K$ (not too large)
Krylov projections (corresponding to $\varphi_p (\tau_k A)w_p$),
performed in smaller Krylov subspaces of dimension $m_1,\cdots,m_K
<m$, which require a total computational cost of
$\mathcal{O}(m^2_1)+\cdots+\mathcal{O}(m^2_K)$. This overall cost is
expected to be less than $\mathcal{O}(m^2)$.  This might be not true
if $K$ is too large, e.g.~the case when the spectrum of $A$ is very
large. This situation, however, is handed by using the adaptive Krylov
algorithm developed in \cite{NW12}, which allows adaptivity of both
the dimension $m$ and the step sizes $\tau_k$.

\subsection{A modification of the $\mathtt{phipm/IOM2}$ routine}
\label{sec:implementation}

Due to the multi-stage structure of exponential Rosenbrock methods, we
have modified the implementation slightly from that used for the
$\epiIII$ method in
\cite{ClancyPudykiewicz13,GaudreaultPudykiewicz16}.  The details of
the previous version (\texttt{phipm/IOM2}) have been previously
described in \cite{GaudreaultPudykiewicz16} (Algorithm 2), which we do
not reproduce here, except for the features relating to our
Rosenbrock-specific modifications.  In particular, each time step of
the $\epiIII$ method requires evaluation of the linear combination
\eqref{eq3.14}; we note that in this formula each of the $\varphi_k$
functions is evaluated at the same argument $A$, and multiplied by a
distinct vector $v_k$.  Moreover, this method utilizes substepping in
time based on formula \eqref{eq3.17}, that effectively evaluates the
$\varphi_k$ functions at scalings of the argument, $\tau_kA$, where
$0< \tau_k\le 1$.  Finally, \texttt{phipm/IOM2} leverages the
recurrence relation \eqref{eq3.18} so that each substep relies on
a set of matrix-vector products, $Aw_{j-1}$.

We now revisit the structure of the multi-stage exponential Rosenbrock
methods from Section \ref{sec:selected_exponential_schemes}.  Unlike
equation \eqref{eq3.14}, these methods compute stages $U_{nj}$ and
time-evolved solutions $u_{n+1}$ by evaluating the $\varphi$ functions
at some scaling of the matrix $A$:
\begin{equation}
\label{eq:exp_multistage}
   w_k = \sum_{l=0}^p \varphi_l(\rho_k\, A) v_l, \quad k=1,\ldots,N_s,
\end{equation}
where the values $\rho_1, \rho_2, \ldots, \rho_{N_s}$ denote an array of ``time''
scaling factors used for each $v_k$ output. For our considered exponential
Rosenbrock schemes, these values are taken from the nodes $c_2, \ldots, c_{s}$.
We further note that the proposed methods only utilize a subset of the $\varphi_l$
functions (e.g., $\mathtt{exprb53}$ uses only $\varphi_1, \varphi_3$
and $\varphi_4$); equivalently, multiple vectors $v_l$ in equation
\eqref{eq:exp_multistage} will be zero.

Due to these structural differences between multi-stage and
multi-step exponential methods, our new \texttt{phipm\_simul\_iom2}
function incorporates two specific changes from the
\texttt{phipm/IOM2} function used in \cite{GaudreaultPudykiewicz16}.

First, we simultaneously compute all of the $w_k$ outputs in equation
\eqref{eq:exp_multistage}, instead of computing these one at a time.
This is accomplished by first requiring that the entire array $\rho_1, \rho_2, \ldots, \rho_{N_s}$  be input to the function.  We then ensure that within the
substepping process \eqref{eq3.17}, these $\rho_j$ values align with a
subset of the internal substep times $t_k$; at each of these moments
the solution vector is stored, and on output the full set
$\{w_k\}_{k=1}^{N_s}$ is returned.  We note that this idea is not new, as
it is similar to an approach used in \cite{TLP12}; ours differs in
that we explicitly stop at each $\rho_k$ instead of computing the $w_k$
vectors through interpolation, thereby guaranteeing no loss of
solution accuracy.

Second, we leverage the fact that some of the input vectors $v_l$ are
identically zero, so within the recursion \eqref{eq3.18} we check
whether $w_{j-1}\ne 0$ before computing the matrix-vector product
$A\,w_{j-1}$.  As matrix-vector products require $\mathcal O(N^2)$
work, while verification that $u\ne 0$ requires only $\mathcal O(N)$,
this can result in significant savings for large $N$, especially when
a significant fraction of these products involve vectors $w_{j-1} = 0$.

Aside from these changes, all components of \texttt{phipm/IOM2} are
directly retained, including the strategy for adapting the substep
size and Krylov subspace size.
We provide pseudocode for this approach in Algorithm~\ref{alg1} below.

\begin{algorithm}[H]
\caption{(\texttt{phipm\_simul\_iom2}): Simultaneously compute all linear combinations in \eqref{eq:exp_multistage}}
\label{alg1}
\begin{list}{$\bullet $}{}
{\small
\item \textbf{Input:}  $A \in \mathbb{R}^{N\times N}$,
$\rho =[\rho_1, \rho_2, \ldots, \rho_{N_s}] $ (where $0<\rho_1 < \ldots <\rho_{N_s}\le 1$),
$v=[v_0,\ldots,v_p]$, $Tol$ (desired tolerance), $m$ (Krylov dimension $\ll N$), and $iom$ (orthogonalization length)
\item \textbf{Initialization:}
$m=1$,\ $m_{max}=100,\  iom=2$,\ $k=0,\ t_k=0,\ y_k=v_0$,\ $\beta=\|v_0\|$,\
$\tau_k=\rho_1$ (initial substep $\tau_0$),\ $\delta =1.4$ (safety factor)
\item \textbf{for $i=1$ to $N_s$ do}\\
\hspace*{0.3cm} $t_{out}=\rho(i)$
\\
\hspace*{0.3cm} \textbf{while} $t_k<t_{out}$
\begin{enumerate}
\item Compute vectors $w_j$ based on \eqref{eq3.18}: $w_0=y_k$, \\
$w_j = \begin{cases} Aw_{j-1}+\sum_{\ell=0}^{p-j} \frac{t^{\ell}_k}{\ell!}v_{j+\ell}, & w_{j-1}\ne 0\\
                                 \sum_{\ell=0}^{p-j} \frac{t^{\ell}_k}{\ell!}v_{j+\ell}, & w_{j-1} = 0\end{cases} \quad j=0,\ldots,p$.
\item Compute $\varphi_p (\tau_k A)w_p$ in \eqref{eq3.17} using the Krylov subspace method:
\begin{enumerate}
  \item [i.] Perform IOM2 method to produce: basis $V_m$ for $K_m$,
  Hessenberg matrix $H_m$, $h_{m+1,m}$, and vector $v_{m+1}$ (as in \cite[Alg.~1]{GaudreaultPudykiewicz16}).
  \item [ii.] Approximate\\
   $\varphi_p (\tau_k A)w_p \approx \beta V_m \varphi_p (\tau_k H_m)e_1 +\beta h_{m+1,m}e^{T}_m \varphi_{p+1} (\tau_k H_m)e_1 v_{m+1}$
 \end{enumerate}
 \item Compute error in the approximation of $\varphi_p (\tau_k A)w_p$:\\
 $\epsilon_m=\beta |h_{m+1,m}| [\varphi_{p+1} (\tau_k H_m)]_{m,1}$
 (as described in  \cite[Sec. 3.2]{NW12})
 \item Update $\tau_{new}$ and $m_{new}$ (as in \cite[Sec. 3.4]{NW12} or \cite[Alg.~3]{GaudreaultPudykiewicz16})
 \item Compute the cost functions $C(\tau_{new},m)$ and $C(\tau, m_{new})$ as in \eqref{eq:cost} and decide whether to update $\tau$ or $m$
 \item \textbf{if} $\frac{t_{out}\|\epsilon_m\|}{\tau_{k} Tol} \leq \delta$:  update $ t_k:=t_k+\tau_k, \ y_k:= y_{k+1}$, $k:=k+1$,
  and\\
  \hspace*{0.3cm} compute $y_{k+1}$ from \eqref{eq3.17}, i.e.,
  $y_{k+1}=\tau^p_k \varphi_p (\tau_k A)w_p+ \sum_{j=0}^{p-i} \frac{\tau^j_k}{j!}w_j$.\\
 \textbf{else}: try again with the revised $\tau$ and $m$ values (back to step 1)
 \end{enumerate}
 \hspace*{0.38cm} \textbf{end}\\
 \hspace*{0.3cm} Store the $i$-th desired linear combination $W(:,i)= \sum_{l=0}^{p} \rho(i)^l  \varphi_l(\rho(i)A)v_l$ \\
 \textbf{end}
 \item \textbf{Output:} $W\in\mathbb{R}^{N \times N_s}$ contains the desired linear
 combination evaluated at the set of times specified by the input array $\rho$.
}
\end{list}
\end{algorithm}
We note that in step 2 (ii), $\varphi_p (\tau_k H_m)e_1$ and
$\varphi_{p+1} (\tau_k H_m)e_1$ are computed simultaneously by computing
$\ee^{\tau_k\hat{H}_m}$, where $\hat{H}_m$ is an augmented matrix of $H_m$.\\
As described in \cite[Sec.~5.3.1]{GaudreaultPudykiewicz16}, the total estimated
cost to advance from the current time step $t_k$ to $t_{out}$ using substeps of
size $\tau$ and a Krylov subspace of size $m$ is
\begin{equation}
\label{eq:cost}
  C(\tau, m)=\nint{\frac{t_{out}-t_k}{\tau}} m(N+n_{A})+2(p-1)(n_{A}+N)+M(m+p+1)^3+(2p+1)N,
\end{equation}
where $n_{A}$ is the number of nonzero entries in $A$ and
$M=\frac{44}{3}+2\nint{\log_{2} \frac{\|\tau_k\hat{H}_m\|}{5.37}}$ is the total
cost of computing $\ee^{\tau_k\hat{H}_m}$.
We further note that the
$\mathtt{phipm\_simul\_iom2}$ function may be downloaded from the Github
repository \url{https://github.com/drreynolds/Phipm_simul_iom}.

\subsection{Details of implementation of $\mathtt{exprb42}$, $\mathtt{pexprb43}$ and $\mathtt{exprb53}$}
We now make use of \texttt{phipm\_simul\_iom2} for implementing the three selected exponential Rosenbrock schemes $\mathtt{exprb42}$, $\mathtt{pexprb43}$ and $\mathtt{exprb53}$.
In the following we will denote $A=\Delta t J_n$ and $v=\Delta t F(u_n)$ for simplicity.\\
\\
\emph{Implementation of $\mathtt{exprb42}$:}
In view of the structure of $\mathtt{exprb42}$ given in  \eqref{eq3.11}, it requires two calls to \texttt{phipm\_simul\_iom2}:
\begin{itemize}
\item[(i)] Evaluate  $y_1=\varphi _{1} (\tfrac{3}{4} A)w_1$ with
  $w_0=0$ and $w_1= \tfrac{3}{4} v$ to get $U_{n2}=u_n+ y_1$.
\item[(ii)] Evaluate $w=\varphi _{1} (A)v_1+ \varphi _{3} (A)v_3$
  with $v_0=v_2=0$, $v_1=v$ and $v_3=\tfrac{32}{9} \Delta t D_{n2}$ to get $u_{n+1}=u_n+ w$.
\end{itemize}
\emph{Implementation of $\mathtt{pexprb43}$:}
Due to the structure of $\mathtt{pexprb43}$ given in  \eqref{eq3.12},
we only need two calls to \texttt{phipm\_simul\_iom2}:
\begin{itemize}
\item[(i)] Evaluate both $y_1=\varphi _{1} (\tfrac{1}{2} A)v$ and
  $z_1=\varphi _{1} (A)v$ simultaneously to compute the stages
  $U_{n2}=u_n+\tfrac{1}{2}y_1$ and $U_{n3}=u_n+z_1$.
\item[(ii)] Evaluate $w=\varphi _{3} (A)v_3+ \varphi _{4} (A)v_4$
  with $v_0=v_1=v_2=0$, $v_3=\Delta t(16D_{n2}-2D_{n3})$ and
  $v_4=\Delta t (-48D_{n2}+12D_{n3})$ to get $u_{n+1}=U_{n3}+w$.
\end{itemize}
\emph{Implementation of $\mathtt{exprb53}$:}
Although $\mathtt{exprb53}$ is a 3-stage scheme like
$\mathtt{pexprb43}$, its structure \eqref{eq3.13} is more
complicated since the $\varphi$-functions in the stage
$U_{n3}$ involve different scaling of $A$. This requires one
additional call to \texttt{phipm\_simul\_iom2} as compared to
$\mathtt{pexprb43}$:
\begin{itemize}
\item[(i)] Evaluate both $y_1=\varphi _{1} (\tfrac{1}{2} A)v$ and
  $z_1=\varphi _{1} (\tfrac{9}{10} A)v$ simultaneously to get $U_{n2}=u_n+ \tfrac{1}{2}y_1$.
\item[(ii)] Evaluate both $y_3=\varphi _{3} (\tfrac{1}{2} A)w_3$ and $z_3=\varphi _{3} (\tfrac{9}{10} A)w_3$ simultaneously with $w_3= hD_{n2}$ to get $U_{n3}=u_n+ \tfrac{9}{10} z_1+ \tfrac{27}{25} y_3+\tfrac{725}{125} z_3$.
\item[(iii)] Evaluate $w=\varphi _{1} (A)v_1+\varphi _{3} (A)v_3+
  \varphi _{4} (A)v_4$  with $v_0=v_2=0$, $v_1=v$,
  $v_3=\Delta t (18D_{n2}-\tfrac{250}{81}D_{n3})$ and $v_4=\Delta t (-60D_{n2}+\tfrac{500}{27}D_{n3})$ to get $u_{n+1}=u_n+w$.
 \end{itemize}
\section{Numerical experiments}
\label{sec:numerical_experiments}

In the following subsections we compare the performance of the
proposed set of exponential integration methods against the $\epiIII$
method previously examined in
\cite{ClancyPudykiewicz13,GaudreaultPudykiewicz16}.  As in those
papers, we focus on a set of standard test problems, originally
introduced by Williamson et al.~\cite{WilliamsonEtAl92}.  These test
problems are ordered in increasing difficulty: the L{\" a}uter test
(with analytical solution) is presented in Section \ref{sec:lauter},
the Rossby-Haurvitz wave test is presented in Section
\ref{sec:rossby_haurvitz}, the isolated mountain test is presented in
Section \ref{sec:mountian}, and the unstable jet test is presented in
Section \ref{sec:jet}.  As with the studies
\cite{ClancyPudykiewicz13,GaudreaultPudykiewicz16}, we performed these
experiments with the shallow water equations on an icosahedral grid.
The state vector for each test problem consisted of four prognostic
variables, defined over grid number 6 with $N_g=$ 40,962 vertices,
amounting to a total of 163,848 unknowns.  All simulations were run in
serial in MATLAB using SMU's \emph{Maneframe2} cluster, where each
node is comprised of dual Intel Xeon E5-2695v4 2.1 GHz 18-core
``Broadwell'' processors with 256 GB of DDR4-2400 memory.

For each test problem, we compare the $\epiIII$, $\exprbIV$, $\pexprbIV$ and
$\exprbV$ methods discussed in Section \ref{sec:exponential_methods}.
We use the \texttt{phipm/IOM2} method for $\epiIII$ and all other methods
use the \texttt{phipm\_simul\_iom2} method, as these proved the most
efficient for each of the respective methods.  Both the
\texttt{phipm/IOM2} and \texttt{phipm\_simul\_iom2} methods were run
using a tolerance of $10^{-4}$, an initial Krylov subspace of
size 1 (this is adapted automatically throughout each run), and an
incomplete orthogonalization length of 2.

The parameter $\gamma_h$ used in the dissipation coefficient
\eqref{eq:niu} in these experiments was set to $0.04\times 10^{-2}$
for all tests except for the unstable jet test in Section
\ref{sec:jet}, where this was increased to $1.25\times 10^{-2}$.

For each test, we present a variety of data to examine the performance
of the proposed exponential integration methods on each problem.  As
both the accuracy and the work required per time step for each method
differ, with higher-order methods generally requiring more work per
step, simple measurements of CPU time or error for the same time step
sizes cannot accurately capture the full ``performance'' picture.
Hence, for each test problem we initially provide convergence plots
showing the error in the height field as a function of time step
size.  Here, we examine only the relative $l_{\infty}$ solution error in the
spatial solution at each time
since results in \cite{GaudreaultPudykiewicz16} showed that these were
essentially equivalent to errors computed in the $l_1$ and $l_2$ norms.
Then using these results, we selected specific time step sizes that
provided similar error thresholds ($10^{-4}$, $10^{-5}$, $10^{-6}$
and $10^{-7}$) for each method on each problem, where in general the
higher-order methods could use larger step sizes than lower-order
methods. Finally, using these method-specific step sizes, we present
fair comparisons of the efficiency (error vs CPU time) of each method
on each test problem.

The question of accuracy is a crucial issue, but it is not always addressed
in the context of meteorological models.  We selected the range $[10^{-7},10^{-4}]$
to reflect the realistic accuracy limits for simulations with atmospheric models
that include dynamics and microphysics (condensation, aerosol and chemistry). The state
vector in such a model may contain elements with a different range of values, and
the accuracy levels used in the article were selected to be sufficient for the
components of the state vector with the smallest values.
We note, however, that even more stringent accuracy thresholds for temporal
integration could be considered when utilizing more ``realistic'' atmospheric
models that include more complicated chemical kinetics and cloud micro-physics.

The most important conclusion from the experiments in this section is
somewhat unsurprising: when lower error is desired, the
higher-order methods are much more efficient than the lower-order
methods.  However, when higher errors are allowed, the lower-order
methods are competitive.  In addition to this conclusion, we find that
even at these larger step sizes the proposed methods exhibit similar
conservation properties, and obtain high-quality solutions, in
comparison with the $\epiIII$ method (and hence with the EPI2,
explicit RK4 and semi-implicit predictor-corrector methods examined in
\cite{ClancyPudykiewicz13,GaudreaultPudykiewicz16,ClancyPudykiewicz13b}).

\subsection{Time dependent zonal flow}
\label{sec:lauter}

\begin{figure}[H]
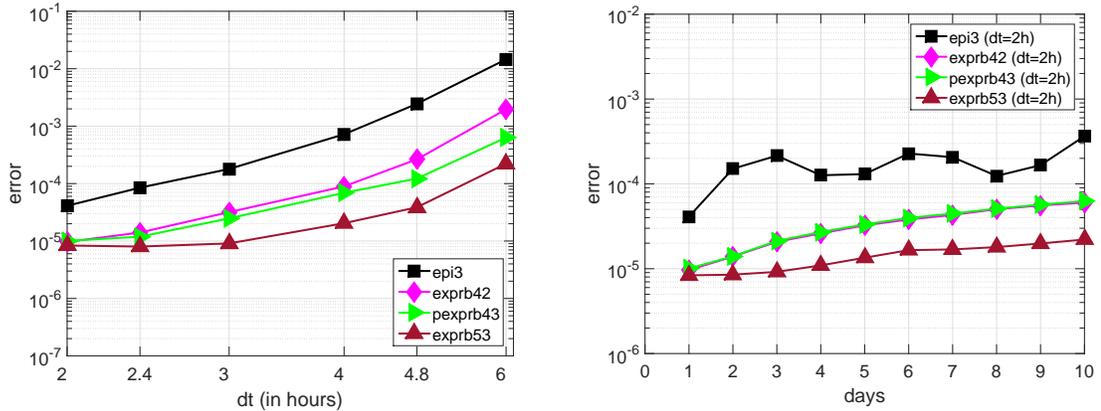

\centering

\begin{tabular}{C{.5\textwidth}C{0pt}C{.5\textwidth}}
{ \includegraphics[width=2.7in, keepaspectratio]{Error_vs_dt_1day_Lauter.eps}}  & &
{ \includegraphics[width=2.7in, keepaspectratio]{Error_vs_dt_upto_10days_Lauter.eps} }
\end{tabular}
\caption{Error measures for the height field in the L\"{a}uter test:
  error after one simulated day using various time step sizes (left),
  and time history of errors over 10 simulated days using a 2 hour
  time-step for each method (right).}
 \label{fig:error_dt_Lauter}
\end{figure}

We begin with a test problem that affords an analytical solution: the
L{\" a}uter test of time dependent zonal flow
\cite{Pudykiewicz2011,LauterHandorfDethloff05}, using the parameters $u_0 = 2\pi a/12$
m/day, $k_1 = 133681$ m$^2$/s$^2$, $k_2 = 10$ m$^2$/s$^2$, and $\alpha
= \pi/4$.
In Figure \ref{fig:error_dt_Lauter} we show two plots.  On the left we
present a standard convergence plot showing error after 1 simulated
day as a function of time step size; on the right we fix each method's
step size to 2 hours, and plot the accumulated solution error over a
10 day simulation.  We note that all curves in the convergence plot
``bottom out'' at an error floor of approximately $10^{-5}$ -- this
corresponds to the spatial approximation error arising from use
  of a single spatial grid in these tests.
We further note that the error history plot shows identical results
for $\epiIII$ as in \cite{GaudreaultPudykiewicz16}, and that the
accumulated error of the proposed methods is uniformly lower than that
of $\epiIII$.

\subsection{Rossby-Haurvitz waves}
\label{sec:rossby_haurvitz}

The shallow water equations used in this study admit several forms
of waves including gravity, Poincar\'e (or inertio-gravity),
Kelvin, and Rossby waves.  In this section we investigate the latter
of these, motivated by their role in our understanding of
atmospheric circulation.  Furthermore, the Rossby wave number four,
used here, is stable as indicated by analytical studies.  The
investigation whether or not our time integration schemes can maintain
this stability is of crucial importance.

As in \cite{GaudreaultPudykiewicz16}, we present a wide range of
results for the Rossby-Haurvitz wave simulation.
These are grouped into two categories: results concerned with solution
error, and results associated with conservation.  When computing
solution error in this problem, we compare against a high-accuracy
reference solution in the absence of an analytical solution.  This
reference solution is computed using the $\pexprbIV$ method with a small
time step size of 30 seconds.
We then compute the $l_{\infty}$ error in the height field by
comparing other solutions against this reference.

\begin{figure}[H]
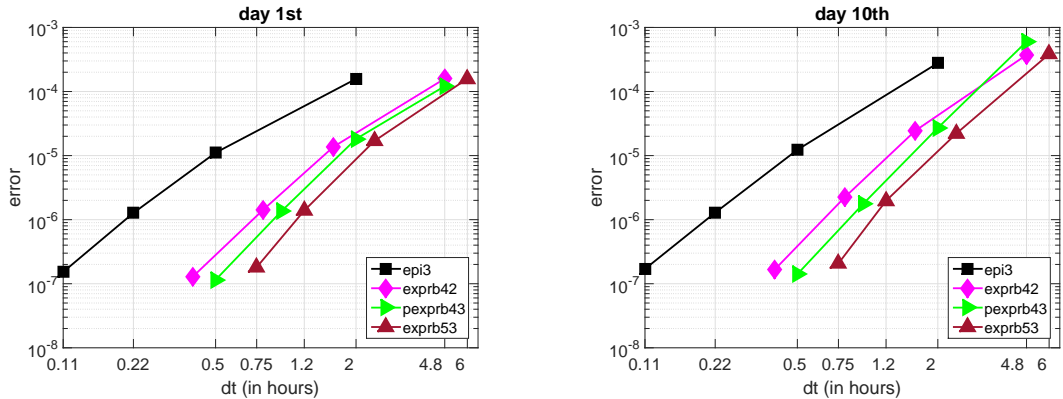

\centering
\begin{tabular}{C{.5\textwidth}R{1pt}C{.5\textwidth}}
{\includegraphics[width=2.5in, keepaspectratio]{Error_vs_dt_1day_Rossby.eps}} & &
{\includegraphics[width=2.5in, keepaspectratio]{Error_vs_dt_10day_Rossby.eps}}
\end{tabular}
\caption{Convergence plots for the Rossby-Haurwitz test.  We plot
  $l_\infty$ error in the height field as a function of time step size
  after 1 simulated day (left) and after 10 simulated days (right).
  The computed order of accuracy for each method is displayed in Table \ref{tab:method_dt_Rossby}.
  }
\label{fig:convergence_Rossby}
\end{figure}

In Figure \ref{fig:convergence_Rossby} we present two ``log-log''
plots of solution error vs time step size for each method; on the left we
compute errors after one simulated day and on the right after 10 simulated
days.  Here, the time step sizes were chosen for each method to
obtain error thresholds approximately equal to $10^{-4}, 10^{-5}, 10^{-6}$
and $10^{-7}$ after one simulated day.  Both plots tell the same story: when
using the same step size the higher-order methods obtain much smaller solution
error, or equivalently, higher-order methods can attain a desired solution
error with much larger time step sizes.

The numerical values represented in the left plot of
Fig.~\ref{fig:convergence_Rossby} are shown in the following Table
\ref{tab:method_dt_Rossby}.  Of particular interest is that the ratio
of time step sizes between $\epiIII$ and the proposed methods increase
precipitously as both the error threshold requirements become more
stringent, and as the method order increases, indicating promise for
higher-order methods when increased accuracy is needed.  The smallest ratio of
2.4 occurs in the $\exprbIV$ and $\pexprbIV$ methods at error threshold
of $10^{-4}$, and the largest ratio of 6.75 occurs in the $\exprbV$
method at an error threshold of $10^{-7}$.

\begin{figure}[H]
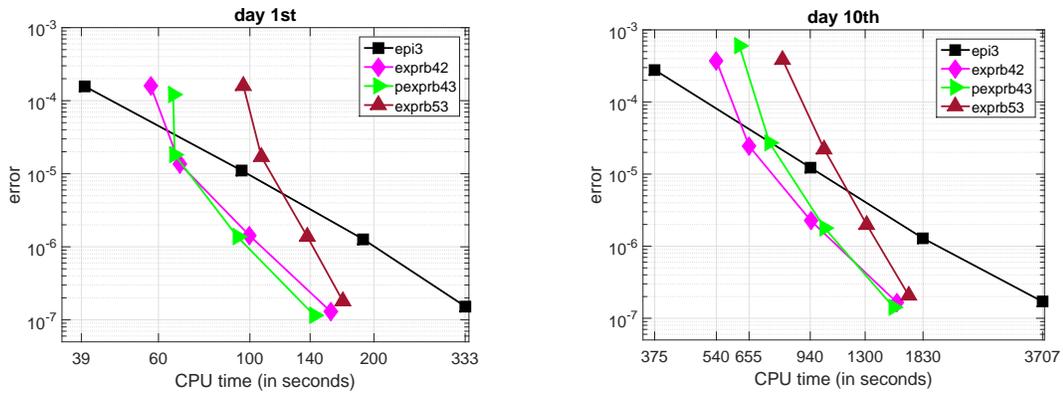

\centering
\begin{tabular}{C{.5\textwidth}R{1pt}C{.5\textwidth}}
{\includegraphics[width=2.5in, keepaspectratio]{Error_vs_CPU_1day_Rossby.eps}} & &
{\includegraphics[width=2.5in, keepaspectratio]{Error_vs_CPU_10day_Rossby.eps}}
\end{tabular}
\caption{Efficiency plots for the Rossby-Haurwitz test.  We plot the
  $l_\infty$ error in the height field as a function of CPU time for a
  variety of time step sizes.  On the left we plot results after one
  simulated day, on the right after 10 simulated days.}
\label{fig:efficiency_Rossby}
\end{figure}

\renewcommand{\arraystretch}{1.2}%
\begin{table}[H]
\centering
\caption{The time step sizes in seconds used for each method to obtain
  desired error thresholds after one simulated day on the Rossby-Haurwitz
  wave test.  In parenthesis, we include the normalized step sizes as
  compared with $\epiIII$ (i.e., we divide by the $\epiIII$ step size
  at each error threshold).
}
\label{tab:method_dt_Rossby}
\small
 \vspace{6pt}
\begin{tabular}{ |c|cc|cc|cc|cc|c| }
\hline
\multirow{2}{*}{Method} &\multicolumn{8}{c|}{Error threshold vs. time step $\Delta t$ }&\multirow{2}{*}{\makecell{Computed\\ order (max.)} } \\
 \cline{2-9}
&\multicolumn{2}{c|}{$10^{-4}$} & \multicolumn{2}{c|}{$10^{-5}$} &\multicolumn{2}{c|}{$10^{-6}$} &\multicolumn{2}{c|}{$10^{-7}$}& \\
 \hline
      $\epiIII$      & 7200   &       & 1800  &         & 800  &        & 400 &                 & 3.06  \\
      $\exprbIV$  & 17280 & (2.4) & 5760 & (3.2) &2880   & (3.6)  & 1440  & (3.6) &  3.46\\
      $\pexprbIV$ & 17280 & (2.4) & 7200   & (4) &3456 & (4.32) & 1800  & (4.5)  &  3.80\\
      $\exprbV$   &21600   & (3)   & 8640 & (4)    & 4320  & (5.4)  & 2700 & (6.75)&  4.34\\
\hline
\end{tabular}
\normalsize
\end{table}
Using these step sizes, in Figure~\ref{fig:efficiency_Rossby} we plot
the efficiency (error vs CPU time) for each method at both one and
ten simulated days.  While qualitatively similar, these plots connote
an expected story: when low-accuracy simulations are desired the
low-order $\epiIII$ method is the most efficient, whereas for error thresholds
of $\sim\!\!10^{-5}$ or tighter the proposed higher-order methods
are more efficient.  We further note that for this test problem, the
increased work-per-step required for the fifth-order $\exprbV$ method
renders it non-competitive for any of the error thresholds examined;
however these results imply that it would become the most efficient
method for error thresholds tighter than $\sim\!\!10^{-8}$.  We note that
this data is also plotted logarithmically, indicating that once the
higher-order methods overtake $\epiIII$, their efficiency benefit over $\epiIII$
improves dramatically.

\begin{figure}[H]
\centering
{\includegraphics[width=5.0in, keepaspectratio]{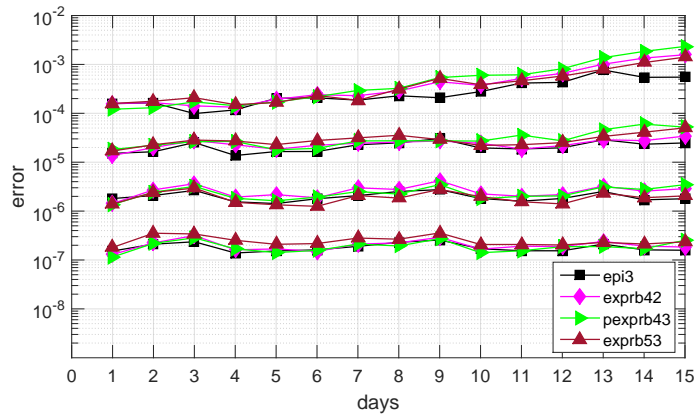}}
\caption{Time history of error for the Rossby-Haurwitz test.  We plot
  the $l_\infty$ error in the height field over 15 simulated days,
  using a variety of time step sizes.
}
\label{fig:error_15days_Rossby}
\end{figure}

Also using the time step sizes from Table~\ref{tab:method_dt_Rossby},
in Figure \ref{fig:error_15days_Rossby} we plot the time history of
the computed error in each method over fifteen simulated days.  While
this plot includes a large amount of data, the salient result is that
for each choice of time step size, all methods exhibit numerical
stability over the full 15 day duration, with errors that do not
accumulate dramatically throughout the simulation.

The following Table \ref{tab:cpu_Rossby} quantifies these results
in slightly more detail.  Specifically, it shows the ratio in CPU
times required for the $\epiIII$ method in comparison with the proposed
methods, for each of the desired error thresholds in performing one
simulated day of the Rossby-Haurwitz wave test.  As described earlier,
for a desired error threshold of $10^{-4}$ the $\epiIII$ method takes less time
than all of the higher-order methods (from 47\% to 69\%); however for
tighter tolerances the higher-order methods are generally faster (with
the exception being $\exprbV$ at $10^{-5}$).  For the tightest
investigated tolerance of $10^{-7}$, each of the proposed methods were
over twice as fast as $\epiIII$, with $\pexprbIV$ attaining the greatest
efficiency improvement.


\renewcommand{\arraystretch}{1.2}%
\begin{table}[H]
\centering
\caption{Speedup factor for the proposed methods in comparison with
  $\epiIII$ for one simulated day of the Rossby-Haurwitz test.  Time
  steps are chosen according to Table \ref{tab:method_dt_Rossby}.}
\label{tab:cpu_Rossby}
\begin{tabular}{l|cccc}
                                  & \multicolumn{4}{c}{Error Threshold}  \\
      Speedup                     & $10^{-4} $  & $10^{-5}$     & $10^{-6}$ &$10^{-7}$     \\ \hline
$\mathtt{cpu_{epi3}/cpu_{exprb42}}$ & 0.69       & 1.41          & 1.89  & 2.12            \\
$\mathtt{cpu_{epi3}/cpu_{pexprb43}}$& 0.61       & 1.45          & 2.02  & 2.32      \\
$\mathtt{cpu_{epi3}/cpu_{exprb53}}$ & 0.47       & 0.89          & 1.37  & 1.98   
\end{tabular}
\end{table}

\begin{figure}[H]
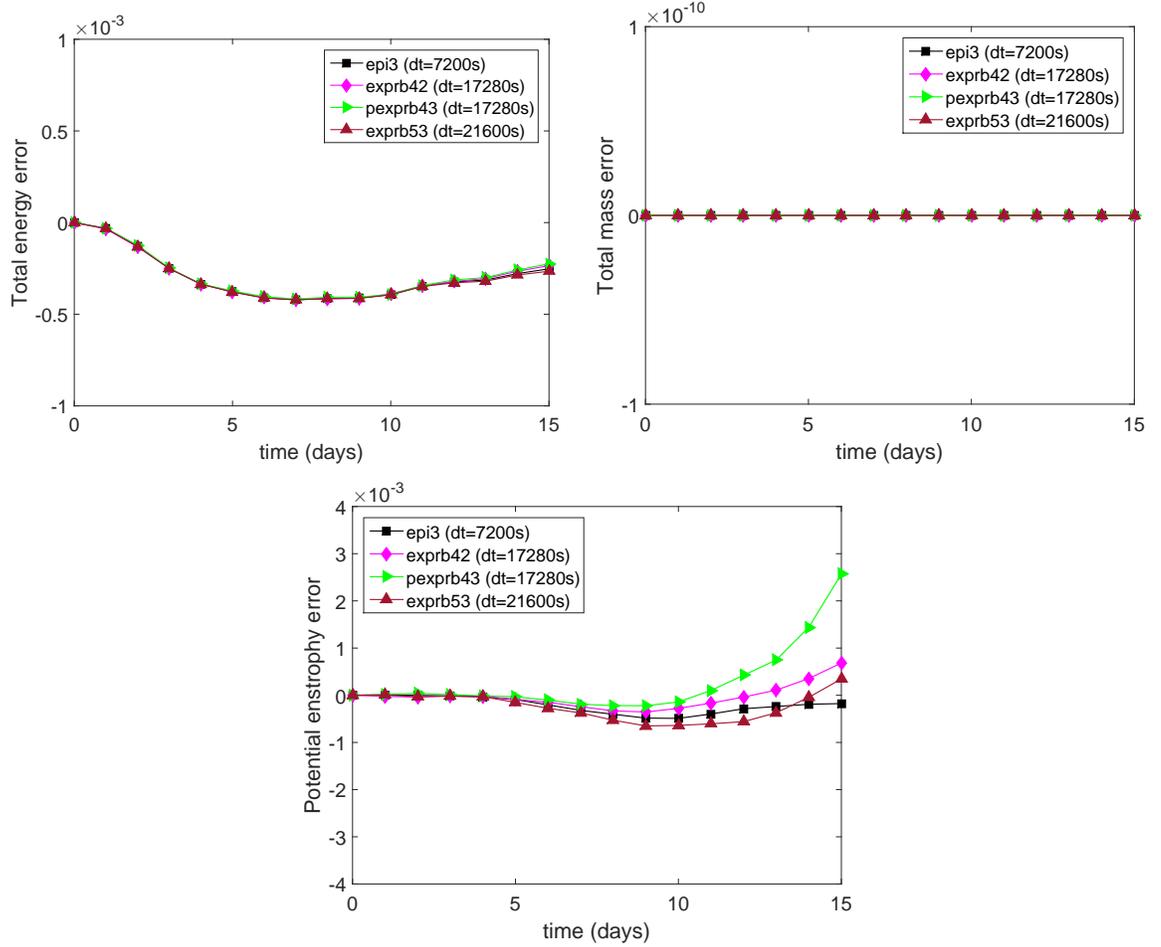

\caption{The normalized conservation errors for total mass, energy and
  potential enstrophy for the Rossby-Haurwitz test.}
\centering
 \label{fig:conservation_epi3_exprb53}
\begin{tabular}{rl}
 {\includegraphics[width=2.9in,  keepaspectratio]{total_energy_error_Rossby.eps}} &
{\includegraphics[width=2.9in, keepaspectratio]{total_mass_error_Rossby.eps}}  \\
 \multicolumn{2}{c}{\includegraphics[width=2.9in, keepaspectratio]{total_enstrophy_error_Rossby.eps}}
\end{tabular}
\end{figure}

We finally consider the conservation properties of these methods on
the Rossyb-Haurwitz test, focusing on total mass, total energy and
potential enstrophy:
\[
  \text{Mass} = h, \quad
  \text{Energy} = gh + \frac{|\bf{u}|^2}{2}, \quad
  \text{Potential Enstrophy} = \frac{\left(\zeta + f\right)^2}{2h}.
\]
In \cite{GaudreaultPudykiewicz16} it was shown that the $\epiIII$
method achieves essentially-perfect conservation of mass, and
conserves both energy and enstrophy to within 0.1\%, in runs over 15
simulated days with time steps of size 2h.  This behavior is
reproduced in Figure \ref{fig:conservation_epi3_exprb53}, which shows
the same behavior for the $\exprbIV$ and $\pexprbIV$ methods using
time steps of size 4.8h, and the $\exprbV$ method using 6h time steps.
We note that these plots show the ``normalized'' conservation
errors, e.g.~the normalized error in conservation of mass is computed
as $\left(\text{Mass}(t) - \text{Mass}(0)\right) /\, \text{Mass}(0)$.

\subsection{The isolated mountain case}
\label{sec:mountian}

We move to the isolated mountain test case from
\cite{WilliamsonEtAl92}, that introduces flows with a complex
vorticity pattern due to perturbation of the flow around an obstacle.
We use an identical setup as in \cite{GaudreaultPudykiewicz16} for
this problem.

\begin{figure}[H]
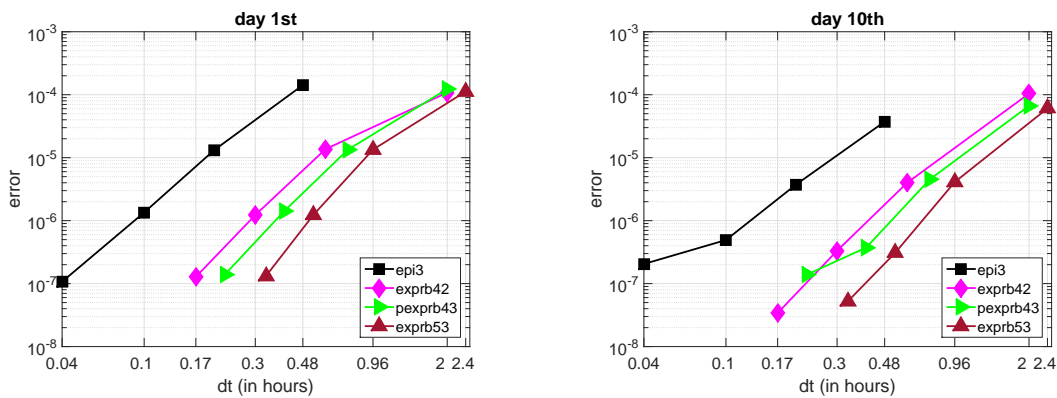

\centering
\begin{tabular}{C{.5\textwidth}R{1pt}C{.5\textwidth}}
{\includegraphics[width=2.5in, keepaspectratio]{Error_vs_dt_1day_Mountain.eps}} & &
{\includegraphics[width=2.5in, keepaspectratio]{Error_vs_dt_10day_Mountain.eps}}
\end{tabular}
\caption{Convergence plots for the mountain test. We plot $\l_\infty$
  error in the height field as a function of time step size after 1
  simulated day (left) and 10 simulated days (right).
  The computed order of accuracy for each method is displayed in Table \ref{tab:method_dt_Mountain}.
  }
 \label{fig:convergence_Mountain}
\end{figure}

\begin{figure}[H]
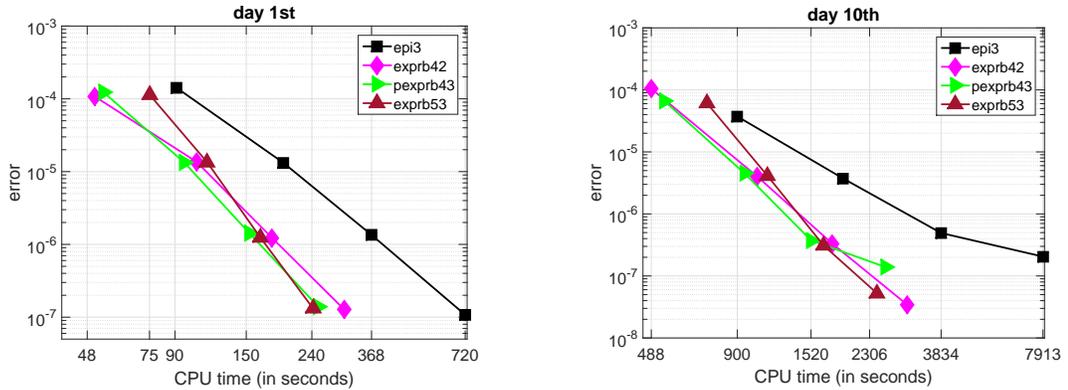

\centering
\begin{tabular}{C{.5\textwidth}R{1pt}C{.5\textwidth}}
{\includegraphics[width=2.5in, keepaspectratio]{Error_vs_CPU_1day_Mountain.eps}} & &
{\includegraphics[width=2.5in, keepaspectratio]{Error_vs_CPU_10day_Mountain.eps}}
\end{tabular}
\caption{Efficiency plots for the mountain test.  We plot the
  $\l_\infty$ error in the height field as a function of CPU time for
  a variety of time step sizes.  Results with one simulated day are on
  the left, and with 10 simulated days are on the right.}
\label{fig:efficiency_Mountain}
\end{figure}

\begin{figure}[H]
\centering
{\includegraphics[width=5in, keepaspectratio]{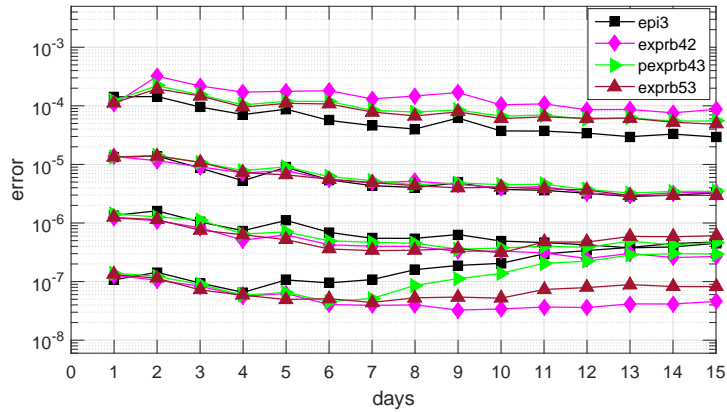}}
\caption{Time history of error for the mountain test.  We plot
  the $l_\infty$ error in the height field over 15 simulated days,
  using a variety of time step sizes.
}
\label{fig:error_15days_Mountain}
\end{figure}

We present an identical set of error-related plots and tables as with
the Rossby-Haurwitz test problem:
\begin{itemize}
\item Figure \ref{fig:convergence_Mountain} plots the convergence of
  each method (similar to Figure \ref{fig:convergence_Rossby}),
\item Table \ref{tab:method_dt_Mountain} shows the step sizes required
  for each method to obtain desired error thresholds after one simulated day
  (similar to Table \ref{tab:method_dt_Rossby}),
\item Figure \ref{fig:efficiency_Mountain} plots the efficiency of
 each method (similar to Figure \ref{fig:efficiency_Rossby}),
\item Figure \ref{fig:error_15days_Mountain} plots the time history of
  the error over 15 simulated days for each method (similar to Figure
  \ref{fig:error_15days_Rossby}),
\item Table \ref{tab:cpu_Mountain} shows the speedup factors for the
  proposed methods in comparison with $\epiIII$ (similar to Table
  \ref{tab:cpu_Rossby}).
\end{itemize}

Again, to compute solution error in this problem we compare against a
high-accuracy reference solution, computed using the $\pexprbIV$ method
with a small time step size of 10 seconds.
We then compute the $l_{\infty}$ error in the height field by
comparing other solutions against this reference.

In Figure \ref{fig:convergence_Mountain} we plot the convergence of
each method on the mountain problam after both 1 and 10 simulated
days.  Two items are particularly interesting.  First, the convergence
curves for the proposed methods are shifted significantly to the right
of the $\epiIII$ curve, indicating that the new methods are much more
accurate for a given time step size on this problem.  Second, in the
10-day results, both the $\epiIII$ and $\pexprbIV$ methods experience
a slight deterioration in convergence at their smallest step sizes.

The numerical values represented in the left plot of
Fig.~\ref{fig:convergence_Mountain} are shown in the following Table
\ref{tab:method_dt_Mountain}.  We note that these time step sizes are
significantly smaller than those for the Rossby-Haurwitz test, due to
the increased difficulty with simulation of flow past an obstacle.
However, similarly to Table \ref{tab:method_dt_Rossby} we note the
trend toward larger step size ratios as error thresholds decrease and as
method order increases.  We also note that here, the smallest time
step ratio between $\epiIII$ and the proposed methods is a factor of
3, and the largest is a factor of 7.5, indicating that the proposed
methods should show even better efficiency improvements for this problem
than the Rossby-Haurwitz test.

\renewcommand{\arraystretch}{1.2}%
\begin{table}[H]
\centering
\caption{The time step sizes in seconds used for each method to obtain
  desired error thresholds after one simulated day on the mountain test.  In
  parenthesis, we include the normalized step sizes as compared with
  $\epiIII$ (i.e.~we divide by the $\epiIII$ step size at each error threshold).
}
\label{tab:method_dt_Mountain}
\vspace{6pt}
\small
\begin{tabular}{ |c|cc|cc|cc|cc|c| }
\hline
\multirow{2}{*}{Method} &\multicolumn{8}{c|}{Error threshold vs. time step $\Delta t$ }&\multirow{2}{*}{\makecell{Computed\\ order (max.)} } \\
 \cline{2-9}
&\multicolumn{2}{c|}{$10^{-4}$} & \multicolumn{2}{c|}{$10^{-5}$} &\multicolumn{2}{c|}{$10^{-6}$} &\multicolumn{2}{c|}{$10^{-7}$}& \\
 \hline
      $\epiIII$   & 1728 &              & 720  &           & 360  &           & 160 &            & 3.12 \\
      $\exprbIV$  &7200  & (4.17) & 2160  & (3)    & 1080   & (3)    & 600  & (3.75) & 3.85  \\
      $\pexprbIV$ & 7200 & (4.17) &2700 & (3.75) &1440   & (4)    & 800  & (4.95)  & 3.95\\
      $\exprbV$   & 8640  & (5)    & 3456 & (4.8)   & 1920 & (5.33) &1200  & (7.5)   & 4.77\\
 \hline
\end{tabular}
\normalsize
\end{table}

Using these step sizes, we again plot the efficiency of each method in
Figure \ref{fig:efficiency_Mountain}.  As anticipated, the proposed
methods are now considerably more efficient than $\epiIII$; 
We again note that the fifth order method is
the most efficient method for error thresholds below $10^{-7}$
both at one and ten simulated days.

Using the time step sizes from Table \ref{tab:method_dt_Mountain}, we
again plot the time history of the computed error in each method over
fifteen simulated days in Figure \ref{fig:error_15days_Mountain}.  We
note that at the smallest error threshold, the artifact noted earlier in
Figure \ref{fig:convergence_Mountain} has more context -- the increase
in error at these runs for the $\epiIII$ and $\pexprbIV$ methods
at the smallest step size begins at day 5, and progressively increases
for the remainder of the run, while the accumulated error for
$\exprbIV$ and $\exprbV$ do not experience a significant increase in
error over the full 15-day simulation.  That said, all methods show
stability at each of these step sizes, even for this more challenging
test.

We again provide in Table \ref{tab:cpu_Mountain} a precise comparison
of the efficiency for each method at these error thresholds,
produced using the results from 10 simulated days above.  As
anticipated from Figure \ref{fig:efficiency_Mountain}, all of the
proposed higher-order methods are more efficient than the $\epiIII$
method, at all error thresholds tested.  This efficiency gap increases as
more accuracy is required, with the higher-order methods uniformly
operating over twice as fast as $\epiIII$ for the error threshold 
of $10^{-6}$, and  $\pexprbIV$ and $\exprbV$ are almost 3 times faster for the tightest
investigated error threshold of $10^{-7}$.

\renewcommand{\arraystretch}{1.2}%
\begin{table}[H]
\centering
\caption{Speedup factor for the proposed methods in comparison with
  $\epiIII$ for ten simulated days of the mountain test.  Time steps
  are chosen according to Table \ref{tab:method_dt_Mountain}.}
\label{tab:cpu_Mountain}
\begin{tabular}{l|cccc}
                                  & \multicolumn{4}{c}{Error Threshold}   \\
      Speedup                      & $10^{-4} $  & $10^{-5}$  & $10^{-6}$ &  $10^{-7}$   \\ \hline
$\mathtt{cpu_{epi3}/cpu_{exprb42}}$  & 1.80       & 1.86       & 2.04      & 2.38   \\
$\mathtt{cpu_{epi3}/cpu_{pexprb43}}$ & 1.70        & 2.04       & 2.40      & 2.87 \\
$\mathtt{cpu_{epi3}/cpu_{exprb53}}$  & 1.21        & 1.72       & 2.22      & 2.98
\end{tabular}
\end{table}

\subsection{The unstable jet case}
\label{sec:jet}

Similarly to \cite{GaudreaultPudykiewicz16}, our last and most
challenging test problem is the unstable jet proposed by Galewsky and
collaborators in 2004 \cite{GalewskyEtAl04}.  This problem generates
vorticity filaments, causing the vorticity gradients to grow
exponentially in time, with a corresponding exponential decay in the
spatial scale required to resolve such filaments.  We counter the
introduction of such small-scale features below the spatial resolution
of the mesh through a slight increase in our dissipation parameter to $\gamma_h=1.25\times 10^{-2}$;
the same value as was used in \cite{ClancyPudykiewicz13}.

It was shown in \cite{GaudreaultPudykiewicz16} that since exponential
integration methods solve the linear portion of the model
analytically, they can successfully resolve the highly-curved flow
field in this test using much longer time steps than would be
otherwise required (e.g.~2 hours for $\epiIII$ versus 30 seconds for
the Robert-Asselin-filtered semi-implicit leapfrog scheme in
\cite{GalewskyEtAl04}).

\begin{figure}[H]
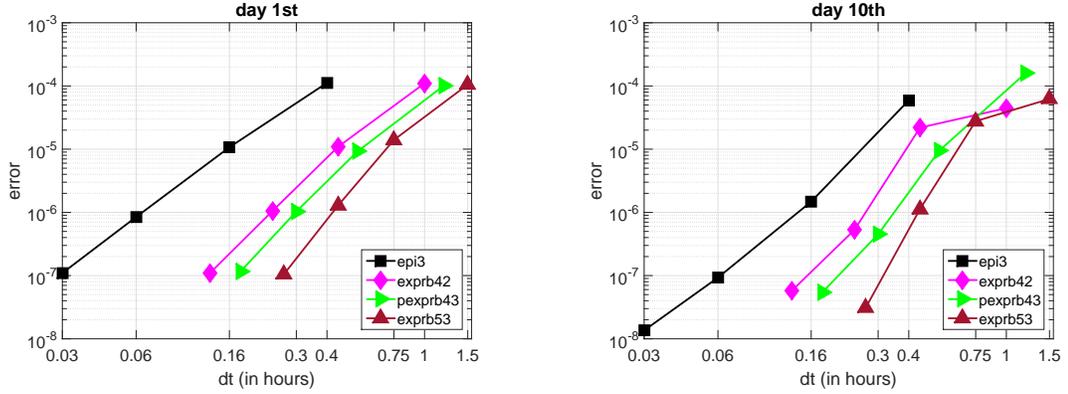

\centering
\begin{tabular}{C{.5\textwidth}R{1pt}C{.5\textwidth}}
 {\includegraphics[width=2.5in, keepaspectratio]{Error_vs_dt_1day_Galewsky.eps}} & &
 {\includegraphics[width=2.5in, keepaspectratio]{Error_vs_dt_10day_Galewsky.eps}}
\end{tabular}
\caption{Convergence plots for the unstable jet test. We plot $\l_\infty$
  error in the height field as a function of time step size after 1
  simulated day (left) and 10 simulated days (right).
  The computed order of accuracy for each method is displayed in Table \ref{tab:method_dt_Galewsky}.
  }
\label{fig:convergence_Galewsky}
\end{figure}

For this problem, we again compare error and efficiency results
for both $\epiIII$ and our proposed higher-order methods,
$\exprbIV$, $\pexprbIV$ and $\exprbV$.  In Figure
\ref{fig:convergence_Galewsky} we show convergence results.  Of note
here is that the error gap between $\epiIII$ is larger after 1
simulated day than after 10 simulated days, although the higher-order
methods can still compute comparably-accurate solutions using
\emph{much} larger time step sizes.

As before, the nuemerical data from the left plot of
Fig~\ref{fig:convergence_Galewsky} are presented in Table
\ref{tab:method_dt_Galewsky}, as well as the
step size ratios when compared to the $\epiIII$ method at the same
error threshold.  Continuing the trend shown in the previous tests, that as
the problems increased in difficulty the higher-order methods promise
enhanced efficiency over $\epiIII$, we note the maximum step size
ratio of 8, again indicating strong potential benefit in using
higher order methods as error threshold requirements tighten.

\renewcommand{\arraystretch}{1.2}%
\begin{table}[H]
\centering
\caption{The time step sizes in seconds used for each method to obtain
  desired error thresholds after one simulated day on the unstable jet test.
  In parenthesis, we include the normalized step sizes as compared with
  $\epiIII$ (i.e.~we divide by the $\epiIII$ step size at each error threshold).
}
\label{tab:method_dt_Galewsky}
\vspace{6pt}
\small
\begin{tabular}{ |c|cc|cc|cc|cc|c| }
\hline
\multirow{2}{*}{Method} &\multicolumn{8}{c|}{Error threshold vs. time step $\Delta t$ }&\multirow{2}{*}{\makecell{Computed\\ order (max.)}} \\
 \cline{2-9}
&\multicolumn{2}{c|}{$10^{-4}$} & \multicolumn{2}{c|}{$10^{-5}$} &\multicolumn{2}{c|}{$10^{-6}$} &\multicolumn{2}{c|}{$10^{-7}$}& \\
\hline
      $\epiIII$      &1440&                & 576 &            & 240  &         & 120&                 &2.94  \\
      $\exprbIV$  & 3600 & (2.5) & 1600   & (2.78)  & 864   & (3.6) & 480 & (4)          &3.85\\
      $\pexprbIV$& 4320 & (3)    & 1920   & (3.33)  & 1080 & (4.5) & 640 & (5.33)     & 4.17\\
      $\exprbV$   & 5400 & (3.75) & 2700 & (4.69)  & 1600 & (6.67)& 960 & (8)         &4.92\\
\hline
\end{tabular}
\normalsize
\end{table}

\begin{figure}[H]
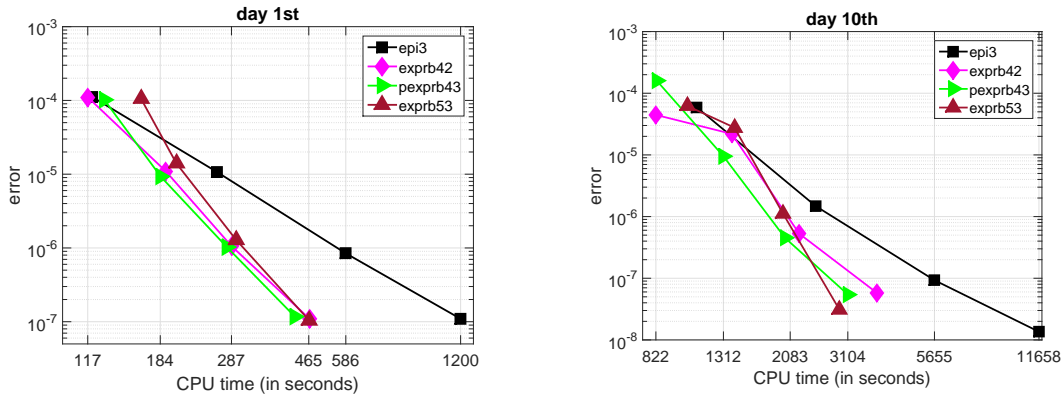

\centering
\begin{tabular}{C{.5\textwidth}R{1pt}C{.5\textwidth}}
{\includegraphics[width=2.5in, keepaspectratio]{Error_vs_CPU_1day_Galewsky.eps}} & &
{\includegraphics[width=2.5in, keepaspectratio]{Error_vs_CPU_10day_Galewsky.eps}}
\end{tabular}
\caption{Efficiency plots for the unstable jet test.  We plot the
  $\l_\infty$ error in the height field as a function of CPU time for
  a variety of time step sizes.  Results with one simulated day are on
  the left, and with 10 simulated days are on the right.}
 \label{fig:efficiency_Galewsky}
\end{figure}

This is again confirmed in the efficiency plots shown in Figure
\ref{fig:efficiency_Galewsky}, where for one simulated day 
both the 4th-order methods 
 are more efficient than $\epiIII$ at all error thresholds, the 5th-order $\exprbV$ beats $\epiIII$ for error thresholds at or below $10^{-5}$.
For ten simulated days the proposed methods beat $\epiIII$ for all
error thresholds at or below $10^{-5}$.  We also note that for the unstable
jet test, the 5th-order $\exprbV$ is competitive with the 4th-order
$\exprbIV$ and $\pexprbIV$ methods for error thresholds at or below
$10^{-6}$ (1 day) and $10^{-5}$ (10 days).

Comparing these CPU times more directly, in Table
\ref{tab:cpu_Galewsky} we provide speedup factors for each method as
compared with $\epiIII$, for each of the investigated error thresholds.
We again note that at error threshold of $10^{-5}$ the proposed methods
are faster than $\epiIII$, with increasing speedups as the desired
solution error decreases, to the point that $\exprbV$ is about 2.5
times faster than $\epiIII$ at the tightest investigated error threshold of
$10^{-7}$.

\renewcommand{\arraystretch}{1.2}%
\begin{table}[H]
\centering
\caption{Speedup factor for the proposed methods in comparison with
  $\epiIII$ for one simulated day of the unstable jet test.
  Time steps are chosen according to Table \ref{tab:method_dt_Galewsky}.}
 \label{tab:cpu_Galewsky}
\begin{tabular}{lcccc}
                                   & \multicolumn{4}{c}{Error Threshold}  \\
      Speedup                                     & $10^{-4}$ & $10^{-5}$ & $10^{-6}$ & $10^{-7}$ \\ \hline
$\mathtt{cpu_{epi3}/cpu_{exprb42}}$  & 1.03     & 1.38      & 2.04     & 2.57      \\
$\mathtt{cpu_{epi3}/cpu_{pexprb43}}$ & 0.93     & 1.42      & 2.10     & 2.81      \\
$\mathtt{cpu_{epi3}/cpu_{exprb53}}$  & 0.73     & 1.29      & 1.98     & 2.57
\end{tabular}
\end{table}

\begin{figure}[H]
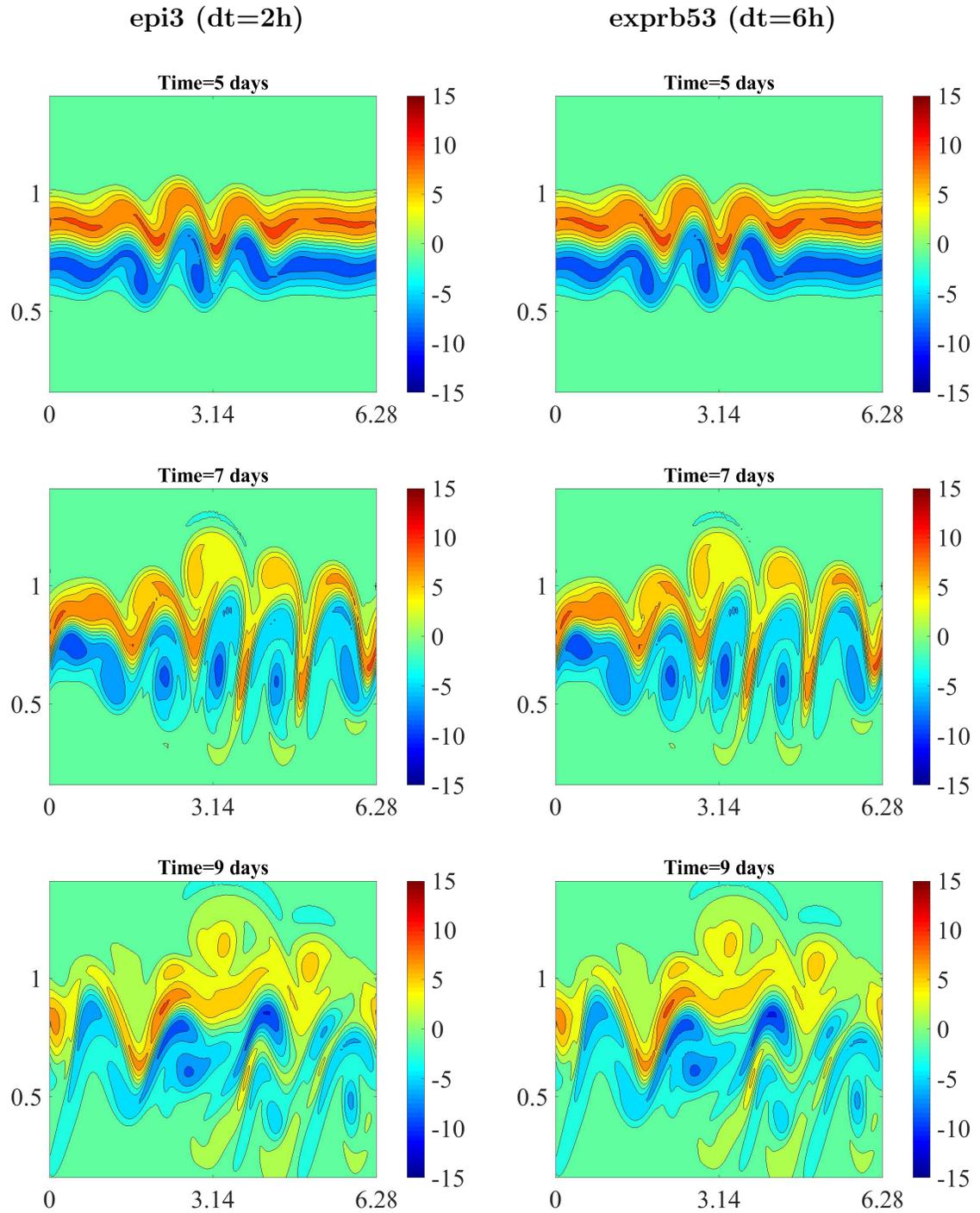

\caption{The vorticity field (scaled by $10^5$) for the unstable jet test.}
 \label{fig:vorti_epi3_exprb53}
\centering
\begin{tabular}{C{.5\textwidth}R{3pt}C{.5\textwidth}}
 \textbf{epi3 (dt=2h)} & & \textbf{exprb53 (dt=6h)}\\
 {\includegraphics[width=2.95in, keepaspectratio]{vorticity_5-epi3.eps}} & &
 {\includegraphics[width=2.95in, keepaspectratio]{vorticity_5-exprb53.eps}} \\[-6pt]
 {\includegraphics[width=2.95in,  keepaspectratio]{vorticity_7-epi3.eps}} & &
 {\includegraphics[width=2.95in, keepaspectratio]{vorticity_7-exprb53.eps}} \\[-6pt]
 {\includegraphics[width=2.95in, keepaspectratio]{vorticity_9-epi3.eps}} & &
 {\includegraphics[width=2.95in, keepaspectratio]{vorticity_9-exprb53.eps}}
\end{tabular}
\end{figure}

We conclude these results by providing contours of the vorticity field
for the unstable jet test in Figure \ref{fig:vorti_epi3_exprb53}.
Here we compare $\epiIII$ against $\exprbV$, where the former uses
time steps of size 2 hours, and the latter of size 6 hours.  These
plots demonstrate qualitative agreement between these two methods
(and in turn the results from \cite{GalewskyEtAl04}).   We note that
the results for $\exprbIV$ and $\pexprbIV$ with time steps of 4h are
essentially identical, and are therefore omitted.

\section{Conclusions and Future Work}
\label{sec:conclusions}

The shallow water equations posed a significant challenge for early
explicit time integration methods, because the CFL criterion resulting
from the stiff gravity waves allowed only very short time steps.
This situation improved with the introduction of semi-implicit
schemes, allowing a sixfold increase in time step size.  Subsequently,
progress was advanced further by early exponential schemes that could
capture the oscillatory term analytically.  The resulting increase in
time step was indeed impressive, allowing a stable, accurate and efficient
integration of the shallow water equations with time steps on the
order of $10^4$ sec on a mesh with horizontal resolution of 100
km (corresponding to a Courant number on the order of 30).
Initial publications in this area addressed the description of the
method in the meteorological context \cite{ClancyPudykiewicz13} and
the efficiency issues \cite{GaudreaultPudykiewicz16}.  In this study,
we continue this effort through formulation of higher order schemes
that approximate the nonlinear part of the problem with increased accuracy.
This property is absolutely crucial in order to make the integrations
with long time steps not just stable but also meaningful.

We provide a detailed summary of both the accuracy and the efficiency
of these methods in Section \ref{sec:numerical_experiments},
indicating that the high order schemes offer a significant advantage
over the methods considered in the initial tests of exponential
schemes in meteorological models, particularly as error
  thresholds are reduced.  Specifically, the higher-order methods were
competitive, or faster than, $\epiIII$ on all test problems when
requiring error thresholds at or below $10^{-5}$.  However, when requiring
an error threshold of $10^{-7}$ the higher-order methods enabled steps from
$3.6$ to $8$ times larger than $\epiIII$ on all test problems.

Furthermore, the efficiency improvements resulting from higher-order
exponential Rosenbrock methods grew markedly stronger as these test
problems increased in nonlinearity.  For the relatively simple flows
exhibited in the L{\" a}uter test in Section \ref{sec:lauter}, the
error threshold needed to be $10^{-5}$ or tighter to benefit from these
higher-order methods, whereas for the turbulent flows exhibited in the
unstable jet test in Section \ref{sec:jet} the higher-order methods
were more efficient than $\epiIII$ at all tested error thresholds.
These results indicate that the increased resolution of the
nonlinear term by the exponential Rosenbrock methods becomes more
critical as flows progress into the nonlinear regime.

One question to be addressed in future research is whether these
benefits can translate to the compressible Euler equations, which form
a fundamental element of current meteorological models.  Preliminary
experiments in this direction are positive.  Even with the second
order EPI2 scheme it was possible to perform an integration of
convective bubbles, reported in \cite {Robert1993}, on a grid with
resolution of 10 m and time step of 5 sec.  The corresponding Courant
number with respect to the acoustic waves was on the order of 150. In
our future work, we will perform systematic tests using our exponential
Rosenbrock methods on this application.  In addition, we will
investigate the role of sources and sinks related to condensation.

We will further explore the benefit of exponential Rosenbrock methods
for atmospheric chemistry, consisting of a large number of coupled
advection--diffusion equations with stiff reaction terms, that still
awaits an ideal time integration algorithm.  The ultimate goal is to
combine, without operator splitting, calculations of the dynamical and
physical processes as well as chemistry in a single stable and
accurate computational framework.
\section*{Acknowledgements}
\label{sec:acknowledgements}

The authors would like to thank the three anonymous referees for their valuable comments and useful suggestions that helped to improve the quality of the paper.

The numerical results were performed on the \emph{Maneframe2} cluster
at Southern Methodist University's Center for Scientific Computation.

\section*{References}
\bibliographystyle{elsarticle-num}
\bibliography{references}

\end{document}